\documentclass[]{interact}

\usepackage{epstopdf}
\usepackage[caption=false]{subfig}

\usepackage[numbers,sort&compress]{natbib}
\bibpunct[, ]{[}{]}{,}{n}{,}{,}
\renewcommand\bibfont{\fontsize{10}{12}\selectfont}

\theoremstyle{plain}
\newtheorem{theorem}{Theorem}[section]
\newtheorem{lemma}[theorem]{Lemma}

\theoremstyle{definition}

\theoremstyle{remark}

\usepackage{tikz,tikz-3dplot}
\usepackage{pgfplots}
\usetikzlibrary{shapes,snakes}
\usepackage{algpseudocode}
\usepackage{algorithm}
\usepackage{url}

\errorcontextlines\maxdimen

\makeatletter
\newcommand*{\algrule}[1][\algorithmicindent]{\makebox[#1][l]{\hspace*{.5em}\thealgruleextra\vrule height \thealgruleheight depth \thealgruledepth}}%
\newcommand*{\thealgruleextra}{}
\newcommand*{\thealgruleheight}{.75\baselineskip}
\newcommand*{\thealgruledepth}{.25\baselineskip}

\newcount\ALG@printindent@tempcnta
\def\ALG@printindent{%
	\ifnum \theALG@nested>0
	\ifx\ALG@text\ALG@x@notext
	\else
	\unskip
	\addvspace{-1pt}
	\ALG@printindent@tempcnta=1
	\loop
	\algrule[\csname ALG@ind@\the\ALG@printindent@tempcnta\endcsname]%
	\advance \ALG@printindent@tempcnta 1
	\ifnum \ALG@printindent@tempcnta<\numexpr\theALG@nested+1\relax
	\repeat
	\fi
	\fi
}%
\usepackage{etoolbox}
\patchcmd{\ALG@doentity}{\noindent\hskip\ALG@tlm}{\ALG@printindent}{}{\errmessage{failed to patch}}
\makeatother

\newbox\statebox
\newcommand{\myState}[1]{%
	\setbox\statebox=\vbox{#1}%
	\edef\thealgruleheight{\dimexpr \the\ht\statebox+1pt\relax}%
	\edef\thealgruledepth{\dimexpr \the\dp\statebox+1pt\relax}%
	\ifdim\thealgruleheight<.75\baselineskip
	\def\thealgruleheight{\dimexpr .75\baselineskip+1pt\relax}%
	\fi
	\ifdim\thealgruledepth<.25\baselineskip
	\def\thealgruledepth{\dimexpr .25\baselineskip+1pt\relax}%
	\fi
	\State #1%
	\def\thealgruleheight{\dimexpr .75\baselineskip+1pt\relax}%
	\def\thealgruledepth{\dimexpr .25\baselineskip+1pt\relax}%
}

\begin{document}


\title{Finding a second Hamiltonian decomposition of a 4-regular multigraph by integer linear programming}

\author{
\name{Andrei V. Nikolaev\thanks{CONTACT A. V. Nikolaev. Email: andrei.v.nikolaev@gmail.com} and Egor V. Klimov}
\affil{P.\,G. Demidov Yaroslavl State University, Yaroslavl, Russia}
}

\maketitle

\begin{abstract}
A Hamiltonian decomposition of a regular graph is a partition of its edge set into Hamiltonian cycles.
We consider the second Hamiltonian decomposition problem: for a 4-regular multigraph find 2 edge-disjoint Hamiltonian cycles different from the given ones.
This problem arises in polyhedral combinatorics as a sufficient condition for non-adjacency in the 1-skeleton of the travelling salesperson polytope.

We introduce two integer linear programming models for the problem based on the classical Dantzig-Fulkerson-Johnson and Miller-Tucker-Zemlin formulations for the travelling salesperson problem.
To enhance the performance on feasible problems, we supplement the algorithm with a variable neighbourhood descent heuristic w.r.t. two neighbourhood structures, and a chain edge fixing procedure.
Based on the computational experiments, the Dantzig-Fulkerson-Johnson formulation showed the best results on directed multigraphs, while on undirected multigraphs, the variable neighbourhood descent heuristic was especially effective.
\end{abstract}


\section{Introduction}

A \textit{Hamiltonian cycle} in the graph $G$ is a cycle that visits each vertex exactly once.
Hamiltonian cycles are named after Sir William Rowan Hamilton, who studied them back in the 1850s and are among the most important and widely studied objects in graph theory.
In particular, the problem of determining whether a graph $G$ contains a Hamiltonian cycle is one of Karp's 21 NP-complete problems \cite{Karp1972}.

A natural generalization of the Hamiltonian cycle problem is the question of finding several edge-disjoint Hamiltonian cycles in the graph $G$.
Of particular interest is the construction of the \textit{Hamiltonian decomposition} of the graph $G$, i.e. partitions of the edge set of $G$ into Hamiltonian cycles.
The classic result by Walecki (1890) states that any complete graph $K_{2n+1}$ with an odd number of vertices has a Hamiltonian decomposition \cite{Alspach2008}.
Later, Walecki's constructions were generalized to other families of graphs, see, for example, the survey \cite{Alspach1990}.
See also the estimate for the number of different Hamiltonian decompositions of a regular graph in \cite{Glebov2017}.
As a combinatorial problem, finding whether a given graph contains a Hamiltonian decomposition or not is NP-complete already for 4-regular undirected and 2-regular directed graphs \cite{Peroche1984}.

In terms of practical applications, the Hamiltonian decomposition is important in the planning of communication networks.
For example, consider the model of all-to-all broadcasting in which each node sends an
identical message to all other nodes in the network.
If the network can be decomposed into edge-disjoint Hamiltonian cycles, then the message traffic will be evenly distributed across all communication links \cite{Rowley1991,Bae2003,Hung2011}.
Among other applications are error-correcting codes \cite{Bailey2009}, privacy preserving distributed data mining \cite{Clifton2002,Dong2010}, and the peripatetic salesperson problem \cite{Krarup1995,DeKort1993}.
Our motivation for this problem comes from the field of polyhedral combinatorics.

\section{Travelling salesperson polytope}

We consider a classic travelling salesperson problem: given a complete weighted graph (or digraph) $K_n=(V,E)$, it is required to find a Hamiltonian cycle of minimum weight.
We denote by $HC_{n}$ the set of all Hamiltonian cycles in $K_{n}$.
With each Hamiltonian cycle $x \in HC_n$ we associate a characteristic vector $\mathbf{x} \in \mathbb{R}^{E}$ by the following rule:
\[
\mathbf{x}_e = 
\begin{cases}
1,& \text{ if the cycle } x \text{ contains an edge } e \in E,\\
0,& \text{ otherwise. }
\end{cases}
\]
The polytope
\[\operatorname{TSP}(n) = \operatorname{conv} \{\mathbf{x} \ | \ x \in HC_n \}\]
is called \textit{the symmetric travelling salesperson polytope}.

\textit{The asymmetric travelling salesperson polytope} $\operatorname{ATSP}(n)$ is defined similarly as the convex hull of characteristic vectors of all possible Hamiltonian cycles in the complete digraph $K_{n}$.

The travelling salesperson polytope was introduced by Dantzig, Fulkerson, and Johnson in their classical work on solving the travelling salesperson problem for 49 US cities by integer linear programming \cite{Dantzig1954}.
State-of-the-art exact algorithms for the travelling salesperson problem are based on a partial description of the facets of the travelling salesperson polytope and the branch and cut method for integer linear programming \cite{Applegate2006}.

The 1-\textit{skeleton} of a polytope $P$ is the graph whose vertex set is the vertex set of $P$ and edge set is the set of geometric edges or one-dimensional faces of $P$.
The study of 1-skeleton is of interest, since, on the one hand, the vertex adjacency can be directly applied to develop simplex-like combinatorial optimization algorithms that move from one feasible solution to another along the edges of the 1-skeleton.
This class includes, for example, the blossom algorithm by Edmonds for constructing maximum matchings \cite{Edmonds1965}, set partitioning algorithm by Balas and Padberg \cite{Balas1975}, Balinski's algorithm for the assignment problem \cite{Balinski1985}, Ikura and Nemhauser's algorithm for the set packing problem \cite{Ikura1985}, etc.
On the other hand, some characteristics of 1-skeleton estimate the time complexity for different computation models and classes of algorithms.
In particular, the diameter (the greatest edge distance between any pair of vertices) is a lower bound for the number of iterations of the simplex-method and similar algorithms \cite{Dantzig1963,Grotchel1985}, while the clique number (the number of vertices in the largest clique) estimates the complexity in the class of \textit{direct type} algorithms based on linear comparisons \cite{Bondarenko1983,Bondarenko2013}.

For some polytopes, there are easily verifiable necessary and sufficient conditions for the vertex adjacency in 1-skeletons. See, for example, vertex covering, partitioning, linear and partial ordering \cite{Hausmann1978}, set covering \cite{Aguilera2017}, constrained assignment problem \cite{Alfakih1998}, fractional stable set polytope \cite{Michini2014}, etc.
Unfortunately, Papadimitriou \cite{Papadimitriou1978} proved that verifying vertex adjacency in 1-skeleton of the travelling salesperson polytope is a hard problem.

\begin{theorem} [Papadimitriou \cite{Papadimitriou1978}]
	The question of whether two vertices of the polytopes $\operatorname{TSP}(n)$ or $\operatorname{ATSP}(n)$ are non-adjacent is NP-complete. 
\end{theorem}

\section{Formulation of the problem}

As a result of the NP-completeness of vertex non-adjacency testing in the 1-skeleton of the travelling salesperson polytope, sufficient conditions for non-adjacency are of particular interest.
In this paper, we consider the most general of the known -- sufficient condition by Rao \cite{Rao1976}.

Let $x = (V,E(x))$ and $y=(V,E(y))$ be two Hamiltonian cycles on the vertex set $V$.
We denote by $x \cup y$ a \textit{union multigraph} $(V,E(x) \cup E(y))$ that contains all edges of both cycles $x$ and $y$.
Note that if two cycles contain the same edge $e$, then both copies of the edge are added to the multigraph $x \cup y$.

\begin{lemma}[Rao \cite{Rao1976}]\label{lemma_sufficient}
	Given two Hamiltonian cycles $x$ and $y$, if the union multigraph $x \cup y$ contains a Hamiltonian decomposition into edge-disjoint cycles $z$ and $w$ different from $x$ and $y$, then the corresponding vertices $\mathbf{x}$ and $\mathbf{y}$ of the polytope $\operatorname{TSP}(n)$ (or $\operatorname{ATSP}(n)$) are not adjacent.
\end{lemma}

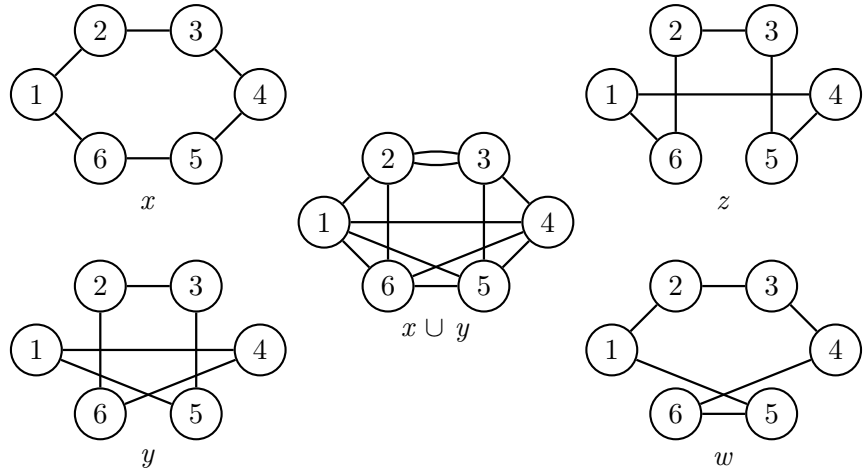
\begin{figure}[t]
	\centering
	\begin{tikzpicture}[scale=0.85]
	\begin{scope}[every node/.style={circle,thick,draw}]
	\node (1) at (0,0) {1};
	\node (2) at (1,1) {2};
	\node (3) at (2.5,1) {3};
	\node (4) at (3.5,0) {4};
	\node (5) at (2.5,-1) {5};
	\node (6) at (1,-1) {6};
	\end{scope}
	\draw [line width=0.3mm] (1) edge (2);
	\draw [line width=0.3mm] (2) edge (3);
	\draw [line width=0.3mm] (3) edge (4);
	\draw [line width=0.3mm] (4) edge (5);
	\draw [line width=0.3mm] (5) edge (6);
	\draw [line width=0.3mm] (6) edge (1);
	\draw (1.75, -1.7) node{\textit{x}};
	
	\begin{scope}[yshift=-4cm]
	\begin{scope}[every node/.style={circle,thick,draw}]
	\node (1) at (0,0) {1};
	\node (2) at (1,1) {2};
	\node (3) at (2.5,1) {3};
	\node (4) at (3.5,0) {4};
	\node (5) at (2.5,-1) {5};
	\node (6) at (1,-1) {6};
	\end{scope}
	\draw [line width=0.3mm] (1) edge (4);
	\draw [line width=0.3mm] (4) edge (6);
	\draw [line width=0.3mm] (6) edge (2);
	\draw [line width=0.3mm] (2) edge (3);
	\draw [line width=0.3mm] (3) edge (5);
	\draw [line width=0.3mm] (5) edge (1);
	\draw (1.75, -1.7) node{\textit{y}};
	\end{scope}
	
	\begin{scope}[xshift=4.5cm,yshift=-2cm]
	\begin{scope}[every node/.style={circle,thick,draw}]
	\node (1) at (0,0) {1};
	\node (2) at (1,1) {2};
	\node (3) at (2.5,1) {3};
	\node (4) at (3.5,0) {4};
	\node (5) at (2.5,-1) {5};
	\node (6) at (1,-1) {6};
	\end{scope}
	\draw [line width=0.3mm] (1) edge (2);
	\draw [line width=0.3mm, bend right=10] (2) edge (3);
	\draw [line width=0.3mm] (3) edge (4);
	\draw [line width=0.3mm] (4) edge (5);
	\draw [line width=0.3mm] (5) edge (6);
	\draw [line width=0.3mm] (6) edge (1);
	\draw [line width=0.3mm] (1) edge (4);
	\draw [line width=0.3mm] (4) edge (6);
	\draw [line width=0.3mm] (6) edge (2);
	\draw [line width=0.3mm, bend left=10] (2) edge (3);
	\draw [line width=0.3mm] (3) edge (5);
	\draw [line width=0.3mm] (5) edge (1);	
	\draw (1.75, -1.7) node{\textit{x $\cup$ y}};
	\end{scope}

	\begin{scope}[xshift=9cm]
	\begin{scope}[every node/.style={circle,thick,draw}]
	\node (1) at (0,0) {1};
	\node (2) at (1,1) {2};
	\node (3) at (2.5,1) {3};
	\node (4) at (3.5,0) {4};
	\node (5) at (2.5,-1) {5};
	\node (6) at (1,-1) {6};
	\end{scope}
	\draw [line width=0.3mm] (1) edge (4);
	\draw [line width=0.3mm] (4) edge (5);
	\draw [line width=0.3mm] (5) edge (3);
	\draw [line width=0.3mm] (3) edge (2);
	\draw [line width=0.3mm] (2) edge (6);
	\draw [line width=0.3mm] (6) edge (1);
	\draw (1.75, -1.7) node{\textit{z}};
	
	\begin{scope}[yshift=-4cm]
	\begin{scope}[every node/.style={circle,thick,draw}]
	\node (1) at (0,0) {1};
	\node (2) at (1,1) {2};
	\node (3) at (2.5,1) {3};
	\node (4) at (3.5,0) {4};
	\node (5) at (2.5,-1) {5};
	\node (6) at (1,-1) {6};
	\end{scope}
	\draw [line width=0.3mm] (1) edge (2);
	\draw [line width=0.3mm] (2) edge (3);
	\draw [line width=0.3mm] (3) edge (4);
	\draw [line width=0.3mm] (4) edge (6);
	\draw [line width=0.3mm] (6) edge (5);
	\draw [line width=0.3mm] (5) edge (1);
	\draw (1.75, -1.7) node{\textit{w}};
	\end{scope}
	\end{scope}
	\end{tikzpicture}
	\caption{The multigraph $x \cup y$ has two different Hamiltonian decompositions}
	\label{image:not_adjacent}
\end{figure}

From a geometric point of view, the sufficient condition means that the segment connecting two vertices $\mathbf{x}$ and $\mathbf{y}$ intersects with the segment connecting two other vertices $\mathbf{z}$ and $\mathbf{w}$ of the polytope $\operatorname{TSP}(n)$ (or $\operatorname{ATSP}(n)$ correspondingly). Hence, the vertices $\mathbf{x}$ and $\mathbf{y}$ are not adjacent. An example of a satisfied sufficient condition is shown in Fig.~\ref{image:not_adjacent}.

We formulate the sufficient condition for vertex non-adjacency of the travelling salesperson polytope as a combinatorial problem.

\vspace{2mm}

\textbf{Second Hamiltonian decomposition problem.}

\textsc{Instance:} let $x$ and $y$ be two Hamiltonian cycles.

\textsc{Question:} does the union multigraph $x \cup y$ contain a pair of edge-disjoint Hamiltonian cycles $z$ and $w$ different from $x$ and $y$?

\vspace{2mm}

Thus, we consider a version of a Hamiltonian decomposition problem of a special form.
By construction, the union multigraph $x \cup y$ always contains the Hamiltonian decomposition into $x$ and $y$.
The question is whether there is a second decomposition into $z$ and $w$.

Unfortunately, checking Rao's sufficient condition is hard.
In particular, Papadimitriou \cite{Papadimitriou1978} proved that the question of whether the union multigraph $x \cup y$ contains a third Hamiltonian cycle $z$ is already an NP-complete problem.
Note that verifying whether a given 4-regular undirected or 2-regular directed multigraph contains a decomposition into edge-disjoint Hamiltonian cycles is also NP-complete \cite{Peroche1984}.

Therefore, instead of Rao's sufficient condition, various polynomially solvable special cases of the vertex non-adjacency problem have been studied in the literature.
In particular, the polynomial sufficient conditions for the pyramidal tours \cite{Bondarenko2018}, pyramidal tours with step-backs \cite{Nikolaev2019}, and pedigrees \cite{Arthanari2006,Arthanari2013} are known.
However, all of them are weaker than the sufficient condition by Rao. 

The second Hamiltonian decomposition problem was introduced in \cite{Kozlova2019} and later studied in \cite{Kostenko2021,Nikolaev2021}:
\begin{itemize}
	\item the simulated annealing algorithm from \cite{Kozlova2019} covers a multigraph with edge-disjoint 2-factors through the reduction to random perfect matching \cite{Tutte1954};
	
	\item the general variable neighbourhood search algorithm from \cite{Nikolaev2021} is a modification of the previous algorithm with two additional neighbourhood structures;
	
	\item the previous version of the algorithm based on integer linear programming is described in the proceedings of the conference ``MOTOR 2021'' \cite{Kostenko2021}.
\end{itemize}

Heuristic algorithms have proven to be very efficient on instances with an existing solution, especially on undirected graphs. However, on instances without a solution, the heuristics face significant difficulties.

In this paper, we consider several exact algorithms for the second Hamiltonian decomposition problem based on the classical ILP-models for the travelling salesperson problem by Dantzig-Fulkerson-Johnson \cite{Dantzig1954} and Miller-Tucker-Zemlin \cite {Miller1960}.
We add to the algorithm the variable neighbourhood descent heuristic adapted from \cite{Nikolaev2021} and enhance it with the chain edge fixing procedure.

Compared to the version in the proceedings of the conference ``MOTOR 2021'' \cite{Kostenko2021}, the local search has been replaced by a variable neighbourhood descent, the Miller-Tucker-Zemlin model has been added for comparison, the implementation of the algorithm has been completely redone and computational experiments were carried out on graphs of different types.

\section{Integer linear programming formulations}

\subsection{Dantzig–Fulkerson–Johnson formulation} \label{Section_DFJ}

Let $x = (V,E(x))$, $y = (V,E(y))$, $x \cup y = (V,E = E(x) \cup E(y))$.
With each edge $e \in E$ we associate the variable
\[
z_e = \begin{cases}
1,& \text{if } e \in z,\\
0,& \text{if } e \in w.
\end{cases}
\]

We adapt the classical ILP-formulation of the travelling salesperson problem by Dantzig, Fulkerson and Johnson \cite{Dantzig1954} into the following ILP-model for the considered second Hamiltonian decomposition problem:
\begin{align}
&\sum_{e \in \delta(v)} z_e = 2, 		&\forall v \in V, \label{DFJ_vertex_degree}\\
&\sum_{e \in E(x) \backslash E(y)} z_e \leq |E(x) \backslash E(y)|-1,  	\label{DFJ_not_x}\\
&\sum_{e \in E(y) \backslash E(x)} z_e \leq |E(y) \backslash E(x)|-1,	\label{DFJ_not_y}\\
&\sum_{e \in E_S} z_e \leq |S|-1, 		&\forall S \subset V,	\label{DFJ_SEC_z}\\
&\sum_{e \in E_S} z_e \geq |E_S| - |S| + 1,		&\forall S \subset V,	\label{DFJ_SEC_w}\\
&z_e \in \{0,1\}, 		&\forall e \in E. 	\label{DFJ_variables}
\end{align}

In the following, we elaborate on the model.
Consider some vertex $v \in V$, we denote by $\delta (v)$ the set of edges incident to $v$.
The \textit{vertex degree constraint} (\ref {DFJ_vertex_degree}) ensures that the degree of each vertex in the cycles $z$ and $w$ is equal to 2.
For directed graphs, the vertex degree constraint is modified so that the indegree and outdegree are equal to 1.

The constraints (\ref{DFJ_not_x})--(\ref{DFJ_not_y}) forbid the Hamiltonian cycles $x$ and $y$ as a solution.
Note that the constraints are imposed only on the unique edges of the cycles since the exchange of multiple edges cannot give a new decomposition.
If we consider a general Hamiltonian decomposition problem of a 4-regular multigraph, then these constraints can be omitted.

Finally, the inequalities (\ref{DFJ_SEC_z})--(\ref{DFJ_SEC_w}) are known as the \textit{subtour elimination constraints} (\textit{SEC}), which forbid solutions consisting of several disconnected tours.
Here $S$ is a subset of $V$, $E_S$ is the set of all edges from $E$ with both vertices belonging to $S$:
\[E_S = \{(u,v) \in E:\ u,v \in S\}.\]

The main problem with the subtour elimination constraints is that there are exponentially many of them: two for each subset of $S \subset V$, i.e. $\Omega (2^{|V|})$.
Therefore, the traditional approach when working with the Dantzig-Fulkerson-Johnson model is to add SEC to the model one at a time as needed.
We start with the relaxed model (\ref{DFJ_vertex_degree})--(\ref{DFJ_not_y}), (\ref{DFJ_variables}) with $O(V)$ constraints, whose integer solutions correspond to a pair of edge-disjoint 2-factors $z$ and $w$.
Find all subtours in $z$ and $w$ and add the constraints (\ref{DFJ_SEC_z})--(\ref {DFJ_SEC_w}) into the model.
We repeat the procedure until the Hamiltonian decomposition is found, or the model is infeasible. This procedure is summarized in Algorithm~\ref{Alg:DFJ}. 

\begin{algorithm}[t]
	\caption{ILP-algorithm based on Dantzig-Fulkerson-Johnson formulation}
	\label{Alg:DFJ}
	\begin{algorithmic}[1]
		\Procedure{DFJ}{$x,y$}
			\State Define current model as (\ref{DFJ_vertex_degree})--(\ref{DFJ_not_y}), (\ref{DFJ_variables})
			\Repeat
				\State $z,w \gets$ an integer point of the current model
				\If{$z$ and $w$ is a Hamiltonian decomposition}
					\State \Return Hamiltonian decomposition $z$ and $w$
				\EndIf
				\State For all subtours in $z$ and $w$ add the SEC (\ref{DFJ_SEC_z})--(\ref{DFJ_SEC_w}) into the model
			\Until the model is infeasible
			\State \Return Hamiltonian decomposition does not exist
		\EndProcedure		
	\end{algorithmic}
\end{algorithm}

\subsection{Miller–Tucker–Zemlin formulation} \label{Section_MTZ}

In 1960, Miller, Tucker, and Zemlin \cite{Miller1960} proposed an alternative integer programming model for the travelling salesperson problem.
The idea was that by introducing additional variables corresponding to the order of vertices in the tour, we can replace the subtour elimination constraints (\ref {DFJ_SEC_z})--(\ref{DFJ_SEC_w}) by equivalent constraints, the number of which is polynomial.
Note that the Miller-Tucker-Zemlin model is designed for the asymmetric travelling salesman problem on directed graphs.

\subsubsection{Directed graphs}

We add to the model (\ref {DFJ_vertex_degree})--(\ref{DFJ_not_y}), (\ref{DFJ_variables}) additional $2(n-1)$  integer variables $\alpha_i, \beta_i$, which denote the order of traversing the vertices in the required tours $z$ and $w$.
Then the subtour elimination constraints (\ref{DFJ_SEC_z})--(\ref{DFJ_SEC_w}) can be replaced with equivalent constraints:
\begin{align}
&2 \leq \alpha_i,\beta_i \leq n,	&\forall i = 2,3,\ldots,n,   \label{MTZ_2_n} \\
&\alpha_i - \alpha_j + n z_e \leq n - 1,	&2 \leq i \neq j \leq n,\ e=(i,j) \in E, \label{MTZ_alpha} \\
&\beta_i - \beta_j + n (1 - z_e) \leq n - 1,	&2 \leq i \neq j \leq n,\ e=(i,j) \in E, \label{MTZ_beta} \\
&\alpha_i, \beta_i \in \mathbb{Z}, 	&\forall i = 2,3,\ldots,n. \label{MTZ_integer}
\end{align}

We assume that both cycles $z$ and $w$ start at vertex number $1$ and then traverse the vertices of the graph in the order of values $\alpha_i$, $\beta_i$.
By constraint (\ref{MTZ_alpha}), if the directed edge $e = (i,j)$ is a part of the tour $z$ (i.e. $z_e = 1$), then
\[\alpha_i \leq \alpha_j - 1,\]
and the vertex $j$ is visited in $z$ after the vertex $i$.
Similarly, the constraint (\ref{MTZ_beta}) guarantees that if $e = (i,j) \in w$ (i.e. $z_e = 0$), then $\beta_i < \beta_j $.

The idea is that if $z$ or $w$ consist of several connected components, then for any numbering $\alpha_i$, $\beta_i$ there is a subtour that does not contain vertex 1.
Such a subtour always contains at least one edge that leads from a higher-numbered vertex to a lower-numbered vertex, which contradicts the constraints (\ref {MTZ_alpha}) and (\ref {MTZ_beta}).

For the directed 2-regular multigraph $x \cup y$, the Miller-Tucker-Zemlin formulation (\ref{DFJ_vertex_degree})--(\ref {DFJ_not_y}), (\ref {DFJ_variables})--(\ref {MTZ_integer}), unlike the Dantzig-Fulkerson-Johnson formulation, contains a linear number of variables and inequalities ($O(V)$).
This allows us to write out all the constraints of the model in an explicit form and find the Hamiltonian decomposition by the ILP-solver in just one iteration, in contrast to Algorithm~\ref{Alg:DFJ}.

\subsubsection{Undirected graphs}

Unfortunately, the Miller-Tucker-Zemlin model cannot simply be transferred to undirected graphs, since it is based on vertex traversal order, i.e. on orientation.

The solution is to duplicate the edges of the graph: for each undirected edge $e = (i,j)$ we create a pair of directed edges $e_1 =(i,j)$ and $e_2 = (j,i)$, only one of which can be visited by the cycles $z$ and $w$.
However, with this approach, the cycles $z$ and $w$ do not cover all the edges of the multigraph $x \cup y$.
Therefore, we can no longer encode each edge of $e \in E$ with just one variable.
With each edge $e = (i,j) \in E$ we associate four variables $z_{i,j}$, $z_{j,i}$, $w_{i,j}$, $w_{j,i}$, such that
\begin{align*}
z_{i,j} &= \begin{cases}
1,& \text{if the tour $z$ traverses the edge in the direction from $i$ to $j$},\\
0,& \text{otherwise};
\end{cases}\\
w_{i,j} &= \begin{cases}
1,& \text{if the tour $w$ traverses the edge in the direction from $i$ to $j$},\\
0,& \text{otherwise}. 
\end{cases}
\end{align*}

The Miller-Tucker-Zemlin formulation for the second Hamiltonian decomposition problem on undirected graphs is obtained by complementing the model (\ref{DFJ_vertex_degree})--(\ref{DFJ_not_y}), (\ref{DFJ_variables})
with additional constraints
\begin{align}
&z_{i,j} + z_{j,i} + w_{i,j} + w_{j,i} = 1,			&\forall (i,j) \in E, \label{MTZ_edge}\\ 
&2 \leq \alpha_i,\beta_i \leq n,	&\forall i = 2,3,\ldots,n,   \label{MTZ_2_n_undirected}\\
&\alpha_i - \alpha_j + n z_{i,j} \leq n - 1,	&2 \leq i \neq j \leq n, \ (i,j) \in E, \label{MTZ_alpha_z}\\
&\beta_i - \beta_j + n w_{i,j} \leq n - 1,	&2 \leq i \neq j \leq n,\ (i,j) \in E, \label{MTZ_beta_w}\\
&z_{i,j}, z_{j,i}, w_{i,j}, w_{j,i} \in \{0,1\}, 		&\forall (i,j) \in E,\\
&\alpha_i, \beta_i \in \mathbb{Z}, 	&\forall i = 2,3,\ldots,n.
\end{align}

The new \textit{edge constraint} (\ref{MTZ_edge}) links the variables $z_{i,j}$ and $w_{i,j}$ and ensures that each edge is included in either $z$ or $w$, and is traversed in only one direction.

The rest of the constraints repeat the model for directed graphs, taking into account that the number of edges and corresponding constraints has increased by 2 times.
Besides, for both variables $z$ and $w$ we double the constraints (\ref{DFJ_not_x})--(\ref{DFJ_not_y}) that forbid the Hamiltonian cycles $x$ and $y$ as a solution, since each cycle can be traversed in two directions: clockwise and counter-clockwise.

The total number of constraints in the Miller-Tucker-Zemlin model for undirected graphs is still linear ($O (V)$). However, there will be about 4 times more variables, and 2.5 times more constraints than in the model for directed graphs.

\subsection{Other formulations}

Note that other ILP-models for the travelling salesperson problem are described in the literature. See, for example, the survey \cite {Orman2007}. However, the Dantzig-Fulkerson-Johnson and Miller-Tucker-Zemlin formulations are the most well-known and applicable in practice.

\section{Heuristics}

To improve the performance, we enhance the iterative Algorithm~\ref{Alg:DFJ} based on Dantzig-Fulkerson-Johnson formulation with the variable neighbourhood descent heuristic.
The neighbourhood structures are modified versions of those used in the GVNS algorithm \cite{Nikolaev2021} with an additional chain edge fixing procedure.

\subsection{Feasible set}

\begin{figure}[h]
	\centering
	\begin{tikzpicture}[scale=1.0]
	\begin{scope}[every node/.style={circle,thick,draw}]
	\node (a1) at (0,0) {1};
	\node (a2) at (1,1) {2};
	\node (a3) at (2.5,1) {3};
	\node (a4) at (3.5,0) {4};
	\node (a5) at (2.5,-1) {5};
	\node (a6) at (1,-1) {6};
	\end{scope}
	
	\node at (1.75, -1.75) {$x \cup y$};
	
	\draw [thick] (a1) edge (a2);
	\draw [thick] (a1) edge (a6);
	\draw [thick] (a2) edge (a6);
	\draw [thick] (a3) edge (a4);
	\draw [thick] (a3) edge (a5);
	\draw [thick] (a4) edge (a5);
	
	\draw [thick, bend left=25] (a1) edge (a2);
	\draw [thick] (a1) edge (a3);
	\draw [thick] (a2) edge (a3);
	\draw [thick] (a4) edge (a6);
	\draw [thick] (a5) edge (a6);
	\draw [thick, bend left=25] (a4) edge (a5);

	\begin{scope}[xshift=-4.75cm]
	\begin{scope}[every node/.style={circle,thick,draw}]
	\node (b1) at (0,0) {1};
	\node (b2) at (1,1) {2};
	\node (b3) at (2.5,1) {3};
	\node (b4) at (3.5,0) {4};
	\node (b5) at (2.5,-1) {5};
	\node (b6) at (1,-1) {6};
	\end{scope}
	
	\draw [thick] (b1) edge (b2);
	\draw [thick] (b1) edge (b6);
	\draw [thick] (b2) edge (b6);
	\draw [thick] (b3) edge (b4);
	\draw [thick] (b3) edge (b5);
	\draw [thick] (b4) edge (b5);
	\node at (1.75, -1.75) {$z$};

	\end{scope}

	\begin{scope}[xshift=4.75cm]
	\begin{scope}[every node/.style={circle,thick,draw}]
	\node (c1) at (0,0) {1};
	\node (c2) at (1,1) {2};
	\node (c3) at (2.5,1) {3};
	\node (c4) at (3.5,0) {4};
	\node (c5) at (2.5,-1) {5};
	\node (c6) at (1,-1) {6};
	\end{scope}
	
	\node at (1.75, -1.75) {$w = (x \cup y) \backslash z$};
	
	\draw [thick, bend left=25] (c1) edge (c2);
	\draw [thick] (c1) edge (c3);
	\draw [thick] (c2) edge (c3);
	\draw [thick] (c4) edge (c6);
	\draw [thick] (c5) edge (c6);
	\draw [thick, bend left=25] (c4) edge (c5);	
	\end{scope}
	\end{tikzpicture}
	\caption{The multigraph $x \cup y$ and its two edge-disjoint 2-factors}
	\label{Fig_2-factors}
\end{figure}
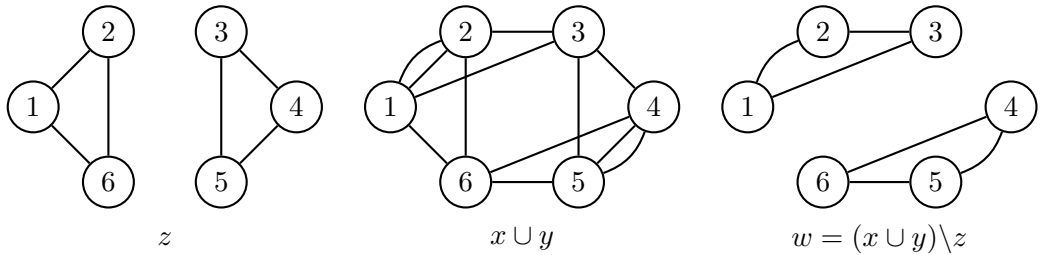

Every solution of the ILP model (\ref{DFJ_vertex_degree})--(\ref{DFJ_not_y}), (\ref{DFJ_variables}) with partial subtour elimination constraints corresponds to the pair $z$ and $w$ of edge-disjoint 2-factors of the multigraph $x \cup y$ (Fig.~\ref{Fig_2-factors}).
Recall that a 2-\textit{factor} is a subgraph of $G$ in which all vertices have degree two.
We compose a set of feasible solutions for the heuristic algorithms from all possible pairs of edge-disjoint 2-factors of the multigraph $x \cup y$.

\subsection{Objective function}

As the objective function, we choose the total number of connected components in the 2-factors $z$ and $w$.
If it equals 2, then $z$ and $w$ are Hamiltonian cycles.

\subsection{Chain edge fixing}

The main difference between the neighbourhood structures in this section and those described in the GVNS algorithm \cite{Nikolaev2021} is the \textit{chain edge fixing} procedure.

We divide the edges of $z$ and $w$ into two classes:
\begin{itemize}
	\item \textit{unfixed} edges that can be moved between $z$ and $w$ to get a neighbouring solution;
	\item edges that are \textit{fixed} in $z$ or $w$ and cannot be moved.
\end{itemize}

\begin{figure}[h]
	\centering
	\begin{tikzpicture}[scale=0.9]
	\begin{scope}[every node/.style={circle,thick,draw}]
	\node (i) at (0,0) {$i$};
	\node (j) at (3,0) {$j$};
	\node (k) at (3,1.5) {$k$};
	\node (h) at (0,-1.5) {$h$};
	\end{scope}
	
	\draw [->,>=stealth,thick] (i) edge node[below]{$z$} (j);
	\draw [->,>=stealth,thick,dashed] (i) edge node[above]{$w$} (k);
	\draw [->,>=stealth,thick,dashed] (h) edge node[below]{$w$} (j);
	
	\end{tikzpicture}
	\caption{Fixing the edge $(i,j)$ in $z$}
	\label{Fig_fixed_edges}
\end{figure}
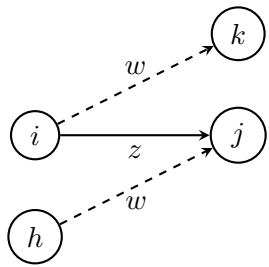

The idea is that one fixed edge starts a recursive chain of fixing other edges.
For example, we consider a directed 2-regular multigraph $x \cup y$ with all indegrees and outdegrees are equal to 2.
Let us choose some edge $(i,j)$ and fix it in the component $z$, then the second edge $(i,k)$ outgoing from $i$ and the second edge $(h,j)$ incoming into $j$ obviously cannot get into $z$. We will fix these edges in $w$ (Fig.~\ref{Fig_fixed_edges}). 
In turn, the edges $(i,k)$ and $(h,j)$, fixed in $w$, start the recursive chains of fixing edges in $z$, etc.

The chain edge fixing procedure for undirected graphs is similar, but with two differences.
First, a recursive chain edge fixing is called when one of the vertices in $z$ or $w$ is incident with two fixed edges.
Second, after chain edge fixing for undirected graphs, the degrees of some vertices in the components $z$ and $w$ may differ from 2.
We call such vertices \textit{broken} and store them in the $brokenList$.
The pseudocode of both procedures is given in Algorithm~\ref{Alg:Chain_Edge_Fixing}.

	\begin{algorithm}[t]
	\caption{Chain edge fixing procedure}\label{Alg:Chain_Edge_Fixing}
	\begin{algorithmic}[1]
		\Procedure{Chain\_Edge\_Fixing\_Directed}{$(i,j)$ in $z$}		
			\State Fix the edge $(i,j)$ in $z$
			\If{the edge $(i,k)$ is not fixed}
				\State \Call{Chain\_Edge\_Fixing\_Directed}{$(i,k)$ in $w$}
			\EndIf
			\If{the edge $(h,j)$ is not fixed}
				\State \Call{Chain\_Edge\_Fixing\_Directed}{$(h,j)$ in $w$}
			\EndIf
		\EndProcedure
		\vspace{2mm}
		\Procedure{Chain\_Edge\_Fixing\_Undirected}{$(i,j)$ in $z$, $brokenList$}		
			\State Update $broken List$ by adding $i,j$ if their degrees are not equal to 2, or excluding if their degrees are restored to 2
			\If{vertex $i$ in $z$ has two incident fixed edges}	
				\State \Call{Chain\_Edge\_Fixing\_Undirected}{unfixed $(i,k)$ in $w$, $brokenList$ }
				\State \Call{Chain\_Edge\_Fixing\_Undirected}{unfixed $(i,p)$ in $w$, $brokenList$ }
			\EndIf
			\If{vertex $j$ in $z$ has two incident fixed edges}
				\State \Call{Chain\_Edge\_Fixing\_Undirected}{unfixed $(j,k)$ in $w$, $brokenList$}
				\State \Call{Chain\_Edge\_Fixing\_Undirected}{unfixed $(j,p)$ in $w$, $brokenList$ }
			\EndIf
		\EndProcedure
	\end{algorithmic}
\end{algorithm}

Note that although chain edge fixing at each step calls several of its recursive copies, the overall complexity of the procedure is linear ($O(V)$), since each edge can be fixed at most once and $|E| = 2|V|$.

\subsection{Local search for directed graphs}

For directed graphs, we supplement the Dantzig-Fulkerson-Johnson formulation with a local search heuristic.

At the preprocessing stage, we fix in $z$ and $w$ a copy of each multiple edges of the union multigraph $x \cup y$, since both copies obviously cannot end up into the same Hamiltonian cycle.

We construct a neighbouring solution by choosing an unfixed edge of $z$, moving it to $w$, and running the chain edge fixing procedure to restore the correct 2-factors $z$ and $w$.
If the number of connected components in $z$ and $w$ has decreased, then we find all subtours in $z$ and $w$, add the corresponding subtour elimination constraints (\ref{DFJ_SEC_z})--(\ref{DFJ_SEC_w}) into the model, proceed to a new solution and continue the local search.
Note that at the beginning we shuffle the edges of $z$ in a random order so that at each run of the local search, the edges are selected with a uniform probability.
This procedure is summarized in Algorithm~\ref{Alg:LS_directed}.

\begin{algorithm}[t]
		\caption{Local search for directed graphs}\label{Alg:LS_directed}
		\begin{algorithmic}[1]
			\Procedure{Local\_Search\_Directed}{$z,w$}
				\State Fix the multiple edges in $z$ and $w$
				\Repeat 
					\State Shuffle the unfixed edges of $z$ in random order	
					\For {each unfixed edge $(i,j)$ in $z$}
						\State \Call {Chain\_Edge\_Fixing\_Directed} {$(i,j)$ in $w$} 	\Comment{Move $(i,j)$ from $z$ to $w$}
						\If {the number of connected components in $z$ and $w$ has decreased}
							\State For all subtours in $z$ and $w$ add the SEC (\ref{DFJ_SEC_z})-(\ref{DFJ_SEC_w}) into the model 
							\State Proceed to a new solution and continue local search
						\EndIf
						\State Restore $z$ and $w$ and unfix all non-multiple edges
					\EndFor
				\Until all edges of $z$ are checked and no improvement found			\Comment{A local minimum}
				\State \Return $z$ and $w$
			\EndProcedure
		\end{algorithmic}
\end{algorithm}

\subsection{Variable neighbourhood descent for undirected graphs}

For undirected graphs, we apply a more complex heuristic with two neighbourhood structures combined in the variable neighbourhood descent approach \cite{Duarte2018}. This is a modified version of the algorithm from \cite {Nikolaev2021}.

The variable neighbourhood search metaheuristic was proposed by Mladenovi\'c and Hansen in 1997 \cite{Mladenovic1997} and has evolved rapidly since then in both its methods and applications. See, for example, surveys by Hansen et al. \cite{Hansen2017, Hansen2019}.

\subsubsection{First neighbourhood structure}

The first neighbourhood for undirected graphs is similar to the neighbourhood structure for directed graphs: we shuffle the edges of $z$ in random order, move one edge from $z$ to $w$, and call the chain edge fixing procedure.

The key difference is that after applying the chain edge fixing procedure on undirected graphs, the components $z$ and $w$ are not necessarily 2-factors. Some broken vertices may remain, the degree of which is not equal to $2$.

In the local search w.r.t. the first neighbourhood structure, in order to restore 2-factors, we extract the next broken vertex from the list, select a random incident unfixed edge, move it between the components $z$ and $w$, and then update the list of broken vertices.
We repeat this operation until the list of broken vertices is empty, after which we get a new feasible solution $z$ and $w$.
If the number of connected components has decreased, we proceed to a new solution. Otherwise, we roll back, restore $z$ and $w$, and try the next edge.
This procedure is summarized in Algorithm~\ref{Alg:VND_first_neighborhood}.

\begin{algorithm}[t]
	\caption{Local search w.r.t. the first neighbourhood structure}\label{Alg:VND_first_neighborhood}
	\begin{algorithmic}[1]
		\Procedure{Local\_Search\_First\_Neighbourhood}{$z,w,attemptLimit$}
			\State Fix the multiple edges in $z$ and $w$
			\Repeat
				\State Shuffle the unfixed edges of $z$ in random order			
				\For {each unfixed edge $(i,j)$ in $z$}
					\State $brokenList \gets \emptyset$
					\State \Call {Chain\_Edge\_Fixing\_Undirected}{$(i,j)$ in $w$, $brokenList$ }
					\For {$i \gets 1$ \textbf{to} $attemptLimit$}				
						\While {$brokenList$ in not emplty}
						\State Extract from $brokenList$ a random broken vertex $i$
						\State Pick a random unfixed edge $(i,k)$ from $w$, if $\deg i = 1$, or from $z$, if $\deg i = 3$
						\State \Call {Chain\_Edge\_Fixing\_Undirected}{$(i,k)$ in $z$ / $w$, $brokenList$ } 
						\EndWhile
						\If {the number of connected components has decreased}
							\State For all subtours in $z$ and $w$ add the SEC (\ref{DFJ_SEC_z})--(\ref{DFJ_SEC_w}) into the model
							\State \Return $z$ and $w$				\Comment{First improvement}
						\EndIf
						\State Restore $z$ and $w$ and unfix all non-multiple edges except $(i,j)$
					\EndFor
					\State Unfix and restore $(i,j)$
				\EndFor
			\Until all edges of $z$ are checked and no improvement found			\Comment{A local minimum}
			\State \Return $z$ and $w$
		\EndProcedure
	\end{algorithmic}
\end{algorithm}

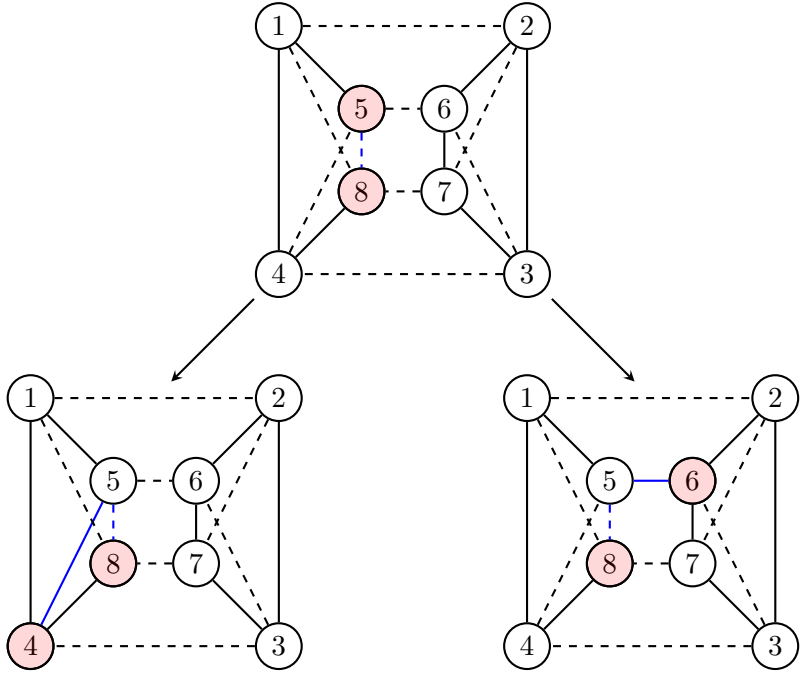
\begin{figure}[t]
	\centering
	\begin{tikzpicture}[scale=1.1]	
	\begin{scope}[every node/.style={circle,draw,thick,inner sep=3pt}]
	\node (1) at (0,3) {$1$};
	\node (2) at (3,3) {$2$};
	\node (3) at (3,0) {$3$};
	\node (4) at (0,0) {$4$};
	\node (5) at (1,2) {$5$};
	\node (6) at (2,2) {$6$};
	\node (7) at (2,1) {$7$};
	\node (8) at (1,1) {$8$};
	\end{scope}
	
	\node [circle,draw,thick,inner sep=3pt,fill=red,fill opacity=0.15] (5) at (1,2) {$5$};
	\node [circle,draw,thick,inner sep=3pt,fill=red,fill opacity=0.15] (8) at (1,1) {$8$};
	
	\draw [thick] (1) -- (5);
	\draw [thick] (8) -- (4);
	\draw [thick] (4) -- (1);
	
	\draw [thick] (2) -- (3);
	\draw [thick] (3) -- (7);
	\draw [thick] (6) -- (7);
	\draw [thick] (2) -- (6);

	\draw [thick,dashed] (1) -- (2);
	\draw [thick,dashed] (2) -- (7);
	\draw [thick,dashed] (7) -- (8);
	\draw [thick,dashed] (1) -- (8);
	
	\draw [thick,dashed] (3) -- (6);
	\draw [thick,dashed] (5) -- (6);
	\draw [thick,dashed] (4) -- (5);
	\draw [thick,dashed] (3) -- (4);
	
	\draw [thick,dashed,blue] (5) -- (8);
	
	\draw [thick,>=stealth,->] (-0.3,-0.3) -- (-1.3,-1.3);
	\draw [thick,>=stealth,->] (3.3,-0.3) -- (4.3,-1.3);

	\begin{scope}[xshift=-3cm,yshift=-4.5cm]
	\begin{scope}[every node/.style={circle,draw,thick,inner sep=3pt}]
	\node (1) at (0,3) {$1$};
	\node (2) at (3,3) {$2$};
	\node (3) at (3,0) {$3$};
	\node (4) at (0,0) {$4$};
	\node (5) at (1,2) {$5$};
	\node (6) at (2,2) {$6$};
	\node (7) at (2,1) {$7$};
	\node (8) at (1,1) {$8$};
	\end{scope}
	
	\node [circle,draw,thick,inner sep=3pt,fill=red,fill opacity=0.15] at (0,0) {$4$};
	\node [circle,draw,thick,inner sep=3pt,fill=red,fill opacity=0.15] at (1,1) {$8$};
	
	\draw [thick] (1) -- (5);
	\draw [thick] (8) -- (4);
	\draw [thick] (4) -- (1);
	
	\draw [thick] (2) -- (3);
	\draw [thick] (3) -- (7);
	\draw [thick] (6) -- (7);
	\draw [thick] (2) -- (6);
	\draw [thick,blue] (4) -- (5);

	\draw [thick,dashed] (1) -- (2);
	\draw [thick,dashed] (2) -- (7);
	\draw [thick,dashed] (1) -- (8);
	
	\draw [thick,dashed] (3) -- (6);
	\draw [thick,dashed] (5) -- (6);
	
	\draw [thick,dashed] (3) -- (4);
	\draw [thick,dashed] (7) -- (8);
	\draw [thick,dashed,blue] (5) -- (8);
	\end{scope}
	
	\begin{scope}[xshift=3cm,yshift=-4.5cm]
	\begin{scope}[every node/.style={circle,draw,thick,inner sep=3pt}]
	\node (1) at (0,3) {$1$};
	\node (2) at (3,3) {$2$};
	\node (3) at (3,0) {$3$};
	\node (4) at (0,0) {$4$};
	\node (5) at (1,2) {$5$};
	\node (6) at (2,2) {$6$};
	\node (7) at (2,1) {$7$};
	\node (8) at (1,1) {$8$};
	\end{scope}
	
	\node [circle,draw,thick,inner sep=3pt,fill=red,fill opacity=0.15] at (2,2) {$6$};
	\node [circle,draw,thick,inner sep=3pt,fill=red,fill opacity=0.15] at (1,1) {$8$};
	
	\draw [thick] (1) -- (5);
	\draw [thick] (8) -- (4);
	\draw [thick] (4) -- (1);
	
	\draw [thick] (2) -- (3);
	\draw [thick] (3) -- (7);
	\draw [thick] (6) -- (7);
	\draw [thick] (2) -- (6);
	\draw [thick,blue] (5) -- (6);

	\draw [thick,dashed] (1) -- (2);
	\draw [thick,dashed] (2) -- (7);
	\draw [thick,dashed] (1) -- (8);
	
	\draw [thick,dashed] (3) -- (6);
	
	\draw [thick,dashed] (4) -- (5);
	\draw [thick,dashed] (3) -- (4);
	\draw [thick,dashed] (7) -- (8);
	\draw [thick,dashed,blue] (5) -- (8);
	\end{scope}
	-\end{tikzpicture}
	\caption{Two ways for restoring a broken vertex $5$ (red) by moving the unfixed (black) incident edges $(5,4)$ or $(5,6)$ from the component $w$ (dashed edges) to the component $z$ (solid edges)
		\label{Fig_restore_broken_vertex}}
\end{figure}

An example of restoring broken vertices is shown in Fig.~\ref{Fig_restore_broken_vertex}. Here the edges of $z$ are solid, the edges of $w$ are dashed, fixed edges are highlighted in blue, and broken vertices are highlighted in red.
The vertex $5$ is broken since its degree in $z$ is $1$.
We can restore the vertex degree in two ways: by moving one of the unfixed incident edges $(5,4)$ or $(5,6)$ from $w$ from $z$.
Note that the third option is theoretically possible: to move the only incident edge $(5,1)$ from $z$ to $w$.
Then the vertex $5$ in $w$ will be incident with 2 fixed edges, and the vertex degree will be automatically restored by the chain edge fixing procedure.
However, in practice, this variant showed worse results, since when restoring one broken vertex we got two new ones, and the list of broken vertices grew rapidly.

Since at each step we choose a random edge to repair the broken vertex, local search w.r.t. to the first neighbourhood for undirected graphs is a randomized algorithm.
Therefore, we run several attempts ($attemptLimit$ parameter) to construct a neighbouring solution.
Thus, the total size of the neighbourhood is $O(V \cdot attemptLimit)$.
While the complexity of constructing one neighbouring solution is $O (V)$, since when repairing broken vertices, at each step at least one edge will be fixed, and the total number of edges is $|E| = 2 |V|$.

\subsubsection{Second neighbourhood structure}

When repairing some broken vertex $i$ at each step, we have several options to move incident edges between the components $z$ and $w$ (Fig.~\ref{Fig_restore_broken_vertex}).
We consider the choice of the edge to move as a branching factor, then all possible choices form a \textit{search tree}, the leaves of which will be either correct 2-factors or spanning graphs that cannot be repaired since some broken vertices contain $3$ or more incident fixed edges.

The local search w.r.t. the first neighbourhood structure explores several random branches in this search tree.
In contrast, a local search w.r.t. the second neighbourhood structure explores all branches of the tree with a depth-first search, but only to a limited depth of recursion.
This technique is called the \textit{bounded search tree} and is the basis of many parametrized algorithms \cite {Cygan2015}.
The pseudocode of the local search w.r.t. the second neighbourhood is given in Algorithm~\ref{Alg:second_neighborhood_structure}.

\begin{algorithm}[t]
		\caption{Local search w.r.t. the second neighbourhood structure}\label{Alg:second_neighborhood_structure}
		\begin{algorithmic}[1]
			\Procedure{Local\_Search\_Second\_Neighbourhood}{$z,w,depthLimit$}
				\State Fix the multiple edges in $z$ and $w$
				\State Shuffle the unfixed edges of $z$ in random order		
				\For {each unfixed edge $(i,j)$ in $z$}
					\State $brokenList \gets \emptyset$
					\State \Call {Chain\_Edge\_Fixing\_Undirected}{$(i,j)$ in $w$,$brokenList$}
					\State \Call{Bounded\_Search\_Tree}{$z,w,depth=1,depthLimit,brokenList$}
					\If {the number of connected components has decreased}
						\State \Return $z$ and $w$		\Comment{First improvement}
					\EndIf
					\State Restore $z$ and $w$ and unfix all non-multiple edges
				\EndFor 
				\State \Return $z$ and $w$			\Comment{A local minimum}
			\EndProcedure
			\vspace{2mm}
			\Procedure{Bounded\_Search\_Tree}{$z,w,depth,depthLimit,brokenList$}
				\If {$z$ and $w$ are 2-factors}
					\If {the number of connected components has decreased}
						\State Proceed to a new solution and exit recursion completely
					\EndIf
					\State \Return				\Comment{Backtrack to the previous level of the search tree}
				\EndIf
				\If {$depth > depthLimit$}
					\State \Return				\Comment{Backtrack to the previous level of the search tree}
				\EndIf
				\State Extract from $brokenList$ a random broken vertex $i$
				\For{each unfixed edge $(i,k)$ incident to $i$ in $z \slash w$}
					\State \Call {Chain\_Edge\_Fixing\_Undirected}{$(i,k)$ in $w \slash z$,$brokenList$}
					\State \Call {Bounded\_Search\_Tree}{$z,w,depth+1,depthLimit,brokenList$}
					\State Restore $z$ and $w$ to the state before the edge $(i,k)$ was moved
				\EndFor
			\EndProcedure
		\end{algorithmic}
\end{algorithm}

The total size of the second neighbourhood is $O (V)$.
The complexity of constructing one neighbouring solution is $O(2^d \cdot V)$, where $d = depthLimit$ is the recursion depth limit, since we are exploring $2^d$ branches of the search tree, and at each branch, the chain edge fixing will be applied to at most all $|E| = 2|V|$ edges.

\subsubsection{Variable neighbourhood descent}
	
For undirected graphs, we combine local search w.r.t. the first and second neighbourhoods in the basic variable neighbourhood descent approach \cite{Duarte2018, Nikolaev2021}: with the first improvement w.r.t. the second neighbourhood, we return to local search w.r.t. the first neighbourhood.
The algorithm stops when either the Hamiltonian decomposition is found or the current feasible solution is a local minimum w.r.t. both neighbourhood structures. This procedure is summarized in Algorithm~\ref{Alg:VND}.

\begin{algorithm}[t]
\caption{Variable neighbourhood descent}\label{Alg:VND}
\begin{algorithmic}[1]
	\Procedure{VND\_Undirected}{$z,w,attemptLimit,depthLimit$}
		\Repeat
			\Repeat
				\State $z$ and $w$ $\gets$ \Call {Local\_Search\_First\_Neighbourhood} {$z,w,attemptLimit$} 
			\Until  $z$ and $w$ is a local minimum w.r.t. the first neighbourhood structure
			\State $z$ and $w$ $\gets$ \Call {Local\_Search\_Second\_Neighbourhood} {$z,w,depthLimit$} 
		\Until  $z$ and $w$ is a local minimum w.r.t. the second neighbourhood structure
		\State \Return $z$ and $w$
	\EndProcedure
\end{algorithmic}
\end{algorithm}

We arrange the neighbourhood structures in order of increasing complexity:
\begin{itemize}
	\item neighbourhood sizes are comparable: $O(V \cdot attemptLimit)$ for the first neighbourhood and $O(V)$ for the second neighbourhood;
	\item the complexity of constructing one neighbouring solution is fundamentally different: $O(V)$ for the first neighbourhood and $O(2^d \cdot V)$ for the second neighbourhood.
\end{itemize}
The order of these neighbourhood structures with computational experiments is investigated in more detail in \cite{Nikolaev2021}.

\subsection{Final algorithm}

We add the above heuristics to Algorithm~\ref{Alg:DFJ} based on the Dantzig-Fulkerson-Johnson formulation between iterations to improve performance on instances with existing Hamiltonian decomposition.
If the ILP-solver returns a pair of edge-disjoint 2-factors $z$ and $w$ that are not a Hamiltonian decomposition, then we call the heuristics to minimize the total number of connected components: local search w.r.t. one neighbourhood structure for directed graphs and variable neighbourhood descent w.r.t. two neighbourhood structures for undirected graphs.

Note that each time the heuristic algorithms improve the solution, we modify the model by adding the corresponding subtour elimination constraints (\ref{DFJ_SEC_z})--(\ref{DFJ_SEC_w}) for all subtours in $z$ and $w$.
Thus, we implement a memory structure that prohibits the algorithm from returning to already investigated feasible solutions.

If the heuristic also fails, we restart the ILP-solver on the modified model and repeat these steps until a Hamiltonian decomposition is found, or the resulting model is infeasible.
This procedure is summarized in Algorithm~\ref{Alg:DFJ+VND}.

\begin{algorithm}[t]
	\caption{ILP-algorithm with the variable neighbourhood descent}
	\label{Alg:DFJ+VND}
	\begin{algorithmic}[1]
		\Procedure{DFJ+VND+Fix}{$x \cup y,attemptLimit,depthLimit$}
			\State Define the current model as (\ref{DFJ_vertex_degree})--(\ref{DFJ_not_y}), (\ref{DFJ_variables})
			\Repeat
				\State $z,w \gets$ an integer point of the current model
				\If{$z$ and $w$ is a Hamiltonian decomposition}
					\State \Return Hamiltonian decomposition $z$ and $w$
				\EndIf
				\State For all subtours in $z$ and $w$ add the SEC (\ref{DFJ_SEC_z})--(\ref{DFJ_SEC_w}) into the model
				\If {the graph is directed}
					\State $z,w \gets$  \Call {Local\_Search\_Directed} {$z,w$};
				\Else
					\State $z,w \gets$  \Call {VND\_Undirected} {$z,w,attemptLimit,depthLimit$};
				\EndIf
				\If{$z$ and $w$ is a Hamiltonian decomposition}
					\State \Return Hamiltonian decomposition $z$ and $w$
				\EndIf
			\Until the model is infeasible
			\State \Return Hamiltonian decomposition does not exist
		\EndProcedure		
	\end{algorithmic}
\end{algorithm}

\section{Computational experiments}

We have chosen 7 algorithms for testing. Three of them are the same for directed and undirected graphs:
\begin{itemize}
	\item MTZ -- Miller-Tucker-Zemlin formulation (Section~\ref{Section_MTZ});
	
	\item DFJ -- Dantzig-Fulkerson-Johnson formulation (Section~\ref{Section_DFJ}) with successive addition of subtour elimination constraints (Algorithm~\ref{Alg:DFJ});
	
	\item SA - simulated annealing algorithm from \cite{Kozlova2019}, which constructs random 2-factors through reduction to perfect matching.
\end{itemize}

On directed graphs, we enhanced the Dantzig-Fulkerson-Johnson model with the local search heuristic:
\begin{itemize}
	\item DFJ+LS -- DFJ formulation with additional local search heuristics and chain edge fixing procedure (Algorithm~\ref{Alg:LS_directed}).
\end{itemize}

On undirected graphs, we enhanced the Dantzig-Fulkerson-Johnson model with three different heuristics:

\begin{itemize}
	\item DFJ+LS-1 -- DFJ formulation with additional local search heuristic w.r.t. the first neighbourhood structure (Algorithm~\ref{Alg:VND_first_neighborhood}), this is the version of the algorithm from the proceedings of the ``MOTOR 2021'' conference \cite {Kostenko2021};
	
	\item DFJ+VND -- DFJ formulation with additional variable neighbourhood descent heuristic w.r.t. two neighbourhood structures (Algorithm~\ref{Alg:DFJ+VND});
	
	\item DFJ+VND+Fix -- the previous algorithm with an additional chain edge fixing procedure (Algorithm~\ref{Alg:Chain_Edge_Fixing}).
\end{itemize}

All the algorithms described in this paper are implemented in .NET Core 3.1, for the SA the existing implementation in Node.js \cite{Nikolaev2019} is taken. Computational experiments were performed on an Intel(R) Core(TM) i7-8750H machine with 2.20GHz CPU and 8GB of RAM. As the ILP-solver we used Gurobi 9.1.2 \cite{Gurobi}.

Note that we did not directly include the general variable neighbourhood search algorithm from \cite{Nikolaev2021} into testing. Its analogue in tests is DFJ+VND, which uses similar neighbourhood structures, but replaces the construction of 2-factors by random matchings with the Dantzig-Fulkerson-Johnson formulation. First, DFJ+VND is an exact algorithm as opposed to the heuristic GVNS. Second, we wanted to minimize the factor of different implementations and programming languages.

We tested the algorithms on $3\,000$ directed and undirected multigraphs from $192$ to $4\,096$ vertices, constructed as unions $x \cup y$ of random Hamiltonian cycles of different nature.

First of all, we need to define a peak in the cycle.
We suppose that the vertices of the graph are labelled from $1$ to $n$. Let $\tau$ be a Hamiltonian cycle.
We denote the successor of $i$-th vertex as $\tau(i)$, and the predecessor as $\tau^{-1}(i)$. A vertex $i$ is called a \textit{peak} if $\tau^{-1}(i)<i$ and $\tau(i)<i$.

\begin{figure}[h]
	\centering
\begin{tikzpicture}[scale=1.0]
\begin{scope}[every node/.style={circle,thick,draw,inner sep=3pt}]
\node (A) at (0,0) {1};
\node (B) at (1,0) {2};
\node (C) at (2,0) {3};
\node (D) at (3,0) {4};
\node (E) at (4,0) {5};
\node (F) at (5,0) {6};
\node (G) at (6,0) {7};
\node (H) at (7,0) {8};
\end{scope}
\draw [thick,->,>=stealth] (A) edge (B);
\draw [thick,->,>=stealth] (B) edge [bend left=50] (D);
\draw [thick,->,>=stealth] (D) edge (E);
\draw [thick,->,>=stealth] (E) edge [bend left=50] (G);
\draw [thick,->,>=stealth] (G) edge (H);
\draw [thick,->,>=stealth] (H) edge [bend left=50] (F);
\draw [thick,->,>=stealth] (F) edge [bend left=50] (C);
\draw [thick,->,>=stealth] (C) edge [bend left=50] (A);
\end{tikzpicture}
	\caption{An example of a pyramidal tour}
\label{image:pyramidal_tour}
\end{figure}
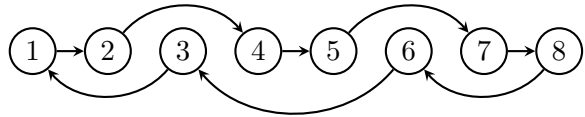

For each graph size, we have constructed $300$ multigraphs $x \cup y$, divided into $3$ groups:
\begin{itemize}
	\item unions of random permutations generated by the Fisher-Yates shuffle algorithm \cite{Knuth1997};
	
	\item unions of \textit{pyramidal tours} with only one peak $n$ (Fig.~\ref{image:pyramidal_tour});
	
	\item unions of Hamiltonian cycles with exactly four peaks.
\end{itemize}

Pyramidal tours are one of the most famous polynomially solvable special cases of the travelling salesperson problem (see, for example, \cite{Gilmore1985}). \textit{Four-peak cycles}, as a generalization of pyramidal tours, are introduced in this paper.

The point is that, on different types of Hamiltonian cycles, the union multigraph $x \cup y$ contains a different number of multiple edges ($|x \cap y|$). 
On the one hand, random permutations generate very few multiple edges: on average $4.04 \pm 2.7$ edges for undirected graphs and $1.99 \pm 1.9$ edges for directed graphs regardless of the graph size. 
On the other hand, pyramidal tours generate a lot of multiple edges: $67.1 \pm 4.4\%$ for undirected graphs and $33.6 \pm 4.1\%$ for directed graphs. Four-peak cycles are an intermediate option: $13.7 \pm 2.1\%$ multiple edges for undirected and $6.8 \pm 1.6\% $ for directed graphs.

Note that the second Hamiltonian decomposition problem on pyramidal tours can be solved in linear time $O (V)$ \cite{Bondarenko2018}. We did not use this fact and chose pyramidal tours as a simple way to generate multigraphs $x \cup y$ with a large number of multiple edges.

For each set of $100$ instances, a limit of 2 hours was set. For the SA algorithm, an additional limit of 60 seconds per test was set, since the heuristic algorithm cannot guarantee that a solution does not exist.
If some algorithm solved less than $25\%$ of the instances from the set, then testing on graphs of this type was stopped.
For each algorithm and each set, the tables with the results indicate how many problems out of $100$ the algorithm solved, the average running time, and the number of iterations with standard deviation.

The results of computational experiments for undirected graphs are given in Table~\ref{table:undirected}. All $3 \,000$ test instances had second Hamiltonian decomposition. Therefore, the table shows the results only for feasible problems.

	\begin{table}[p]
		\centering
		\caption{Computational results for undirected graphs}
		\label{table:undirected}
		\resizebox{\textwidth}{!}{%
			\begin{tabular}{|*{11}{r|}}
				\hline
				& & \multicolumn{3}{c|}{MTZ} & \multicolumn{3}{c|}{DFJ} & \multicolumn{3}{c|}{SA}\\ \hline
				$|V|$ & $|x \cap y|$ & N & time (s) & Iter & N & time (s) & Iter & N & time (s) & Iter  \\ 
				\hline
				
				& $4.4 \pm 2.9$ & $100$ & $5.43 \pm 4.82$ & $1$ & $100$ & $0.91\pm 2.61$ & $27.6 \pm 27.4$ & $100$ & $0.99 \pm 1.69$ & $116 \pm 187$\\
				192 & $26.3 \pm 6.7$ & $99$ & $22.5 \pm 24.1$ & $1$ & $100$ & $48.27 \pm 67.4$ & $135.5 \pm 128.5$ & $19$ & $9.75 \pm 5.18$ & $1275 \pm 716$\\
				& $132 \pm 15.8$ & $1$ & $654.5 \pm 0.00$ & $1$ & $100$ & $0.15 \pm 0.06$ & $15.8 \pm 3.8$ & $0$ & $-$ & $-$\\
				
				\hline
				& $3.8 \pm 2.5$ & $100$ & $16.4 \pm 13.5$ & $1$ & $100$ & $6.71 \pm 29.4$ & $30.9 \pm 35.5$ & $100$ & $2.14 \pm 2.27$ & $137 \pm 135$\\
				256 & $36.5 \pm 8.0$ & $52$ & $33.5 \pm 34.8$ & $1$ & $20$ & $358.6 \pm 565.1$ & $336.6 \pm 352.0$ & $-$ & $-$ & $-$\\
				& $173 \pm 15.2$ & $-$ & $-$ & $-$ & $100$ & $0.32 \pm 0.13$ & $20.1 \pm 4.8$ & $-$ & $-$ & $-$\\
				\hline

				& $3.9 \pm 2.7$ & $62$ & $118.8 \pm 70.3$ & $1$ & $100$ & $7.39 \pm 17.3$ & $41.8\pm 40.2$ & $99$ & $6.62 \pm 6.58$ & $202 \pm 185$\\
				384 & $52.9 \pm 10.5$ & $10$ & $186.1 \pm 99.0$ & $1$ & $0$ & $-$ & $-$ & $-$ & $-$ & $-$\\
				& $261 \pm 18.9$ &$-$ & $-$ & $-$ & $100$ & $0.90 \pm 0.35$ & $27.5 \pm 6.3$ & $-$ & $-$ & $-$\\
				\hline

				& $4.0 \pm 2.7$ & $14$ & $465.3 \pm 222.5$ & $1$ & $100$ & $54.30 \pm 137.3$ & $61.6 \pm 67.2$ & $93$ & $14.6 \pm 13.5$ & $261 \pm 226$\\
				512 & $72.2 \pm 11.4$ & $-$ & $-$ & $-$ & $-$ & $-$ & $-$ & $-$ & $-$ & $-$\\
				& $342 \pm 24.5$ & $-$ & $-$ & $-$ & $100$ & $2.03 \pm 0.74$ & $34.9 \pm 7.5$ & $-$ & $-$ & $-$\\
				\hline

				& $3.9 \pm 2.6$ & $-$ & $-$ & $-$ & $76$ & $87.59 \pm 205.1$ & $60.4 \pm 66.3$ & $77$ & $24.8 \pm 18.8$ & $220 \pm 154$\\
				768 & $105 \pm 13.2$ & $-$ & $-$ & $-$ & $-$ & $-$ & $-$ & $-$ & $-$ & $-$\\
				& $516 \pm 29.2$ & $-$ & $-$ & $-$ & $100$ & $7.16 \pm 2.48$ & $50.3 \pm 9.6$ & $-$ & $-$ & $-$\\
				\hline

				& $4.0 \pm 2.7$ & $-$ & $-$ & $-$ & $16$ & $435.9 \pm 758.4$ & $101.4 \pm 111.5$ & $51$ & $24.6 \pm 17.0$ & $141 \pm 89.6$\\
				1024 & $137 \pm 16.6$ & $-$ & $-$ & $-$ & $-$ & $-$ & $-$ & $-$ & $-$ & $-$\\
				& $681 \pm 35.8$ & $-$ & $-$ & $-$ & $100$ & $20.13 \pm	5.57$ & $69.17 \pm 9.9$ & $-$ & $-$ & $-$\\
				\hline

				& $4.6 \pm 3.1$ & $-$ & $-$ & $-$ & $-$ & $-$ & $-$ & $20$ & $33.5 \pm 17.6$ & $90.8 \pm 26.6$\\
				1536 & $212 \pm 18.8$ & $-$ & $-$ & $-$ & $-$ & $-$ & $-$ & $-$ & $-$ & $-$\\
				& $1017 \pm 39.9$ & $-$ & $-$ & $-$ & $100$ & $70.34 \pm 15.5$ & $100.9 \pm 15.5$ & $-$ & $-$ & $-$\\
				\hline

				& $3.7 \pm 2.6$ & $-$ & $-$ & $-$ & $-$ & $-$ & $-$ & $-$ & $-$ & $-$\\
				2048 & $284 \pm 22.4$ & $-$ & $-$ & $-$ & $-$ & $-$ & $-$ & $-$ & $-$ & $-$\\
				& $1366 \pm 48.6$ & $-$ & $-$ & $-$ & $41$ & $174.3 \pm 39.5$ & $130.0 \pm 21.5$ & $-$ & $-$ & $-$\\
				\hline

				& $3.7 \pm 2.5$ & $-$ & $-$ & $-$ & $-$ & $-$ & $-$ & $-$ & $-$ & $-$\\
				3072 & $410 \pm 28.8$ & $-$ & $-$ & $-$ & $-$ & $-$ & $-$ & $-$ & $-$ & $-$\\
				& $2036 \pm 46.1$ & $-$ & $-$ & $-$ & $7$ & $935.1 \pm 76.9$ & $205.5 \pm 4.8$ & $-$ & $-$ & $-$\\
				\hline

				& $4.3 \pm 2.8$ & $-$ & $-$ & $-$ & $-$ & $-$ & $-$ & $-$ & $-$ & $-$\\
				4096 & $553 \pm 33.6$ & $-$ & $-$ & $-$ & $-$ & $-$ & $-$ & $-$ & $-$ & $-$\\
				& $2742 \pm 67.0$ & $-$ & $-$ & $-$ & $-$ & $-$ & $-$ & $-$ & $-$ & $-$\\
				\hline
				\hline

				& & \multicolumn{3}{c|}{DFJ + LS-1} & \multicolumn{3}{c|}{DFJ+VND} & \multicolumn{3}{c|}{DFJ+VND+Fix} \\ \hline
				$|V|$ & $|x \cap y|$ & N & time (s) & Iter & N & time (s) & Iter & N & time (s) & Iter\\ 
				\hline
				& $4.4 \pm 2.9$ & $100$ & $0.01\pm 0.01$ & $1.04 \pm  0.2$ & $100$ & $0.01\pm 0.01$ & $1.03 \pm 0.0$ & $100$ & $0.01\pm 0.01$ & $1.01 \pm 0.0$\\
				192 & $26.3 \pm 6.7$ & $100$ & $0.03 \pm 0.02$ & $1.33 \pm 0.6$ & $100$ & $0.03 \pm 0.03$ & $1.25 \pm 0.5$ & $100$ & $0.03 \pm 0.03$ & $1.16 \pm 0.4$\\
				& $132 \pm 15.8$ & $100$ & $0.09 \pm 0.09$ & $2.41 \pm 1.6$ & $100$ & $0.04 \pm 0.02$ & $1.06 \pm 0.3$ & $100$ & $0.04 \pm 0.02$ & $1$\\
				\hline
				
				\hline
				& $3.8 \pm 2.5$ & $100$ & $0.02\pm 0.02$ & $1.07\pm 0.3$ & $100$ & $0.02\pm 0.02$ & $1.04\pm 0.2$ & $100$ & $0.02\pm 0.02$ & $1.07 \pm 0.2$\\
				256 & $36.5 \pm 8.0$ & $100$ & $0.05 \pm 0.06$ & $1.43 \pm 0.8$ & $100$ & $0.04 \pm 0.05$ & $1.18 \pm 0.5$ & $100$ & $0.05 \pm 0.06 $ & $1.18 \pm 0.4$\\
				& $173 \pm 15.2$ & $100$ & $0.19 \pm 0.17$ & $2.69 \pm 1.8$ & $100$ & $0.08 \pm 0.05$ & $1.13 \pm 0.4$ & $100$ & $0.07 \pm 0.02$ & $1$\\
				\hline

				& $3.9 \pm 2.7$ & $100$ & $0.04\pm 0.04$ & $1.10 \pm 0.3$ & $100$ & $0.03\pm 0.03$ & $1.05 \pm 0.2$ & $100$ & $0.03\pm 0.03$ & $1.03 \pm 0.2$\\
				384 & $52.9 \pm 10.5$ & $100$ & $0.20 \pm 0.26$ & $1.91 \pm 1.3$ & $100$ & $0.17 \pm 0.20$ & $1.44 \pm 0.8$ & $100$ & $0.14 \pm 0.14$ & $1.32 \pm 0.5$\\
				& $261 \pm 18.9$ & $100$ & $0.60 \pm 0.49$ & $3.34 \pm 1.9$ & $100$ & $0.24 \pm 0.18$ & $1.25 \pm 0.6$ & $100$ & $0.21 \pm 0.08$ & $1$\\
				\hline

				& $4.0 \pm 2.7$ & $100$ & $0.07\pm 0.06$ & $1.12 \pm 0.4$ & $100$ & $0.06 \pm 0.05$ & $1.02 \pm 0.2$ & $100$ & $0.06 \pm 0.07$ & $1.03 \pm 0.2$\\
				512 & $72.2 \pm 11.4$ & $100$ & $0.33 \pm 0.39$ & $1.92 \pm 1.1$ & $100$ & $0.25 \pm 0.30$ & $1.28 \pm 0.6$ & $100$ & $0.29 \pm 0.30$ & $1.34 \pm 0.6$\\
				& $342 \pm 24.5$ & $100$ & $1.68 \pm 1.74$ & $4.14 \pm 2.7$ & $100$ & $0.39 \pm 0.13$ & $1.07 \pm 0.3$ & $100$ & $0.44 \pm 0.14$ & $1$\\
				\hline
				
				& $3.9 \pm 2.6$ & $100$ & $0.14 \pm 0.17$ & $1.10 \pm 0.3$ & $100$ & $0.11 \pm 0.11$ & $1.03 \pm 0.2$ & $100$ & $0.13 \pm 0.14$ & $1.03 \pm 0.2$\\
				768 & $105 \pm 13.2$ & $100$ & $1.04 \pm 1.49$ & $1.97 \pm 1.4$ & $100$ & $0.93 \pm 1.15$ & $1.56 \pm 0.9$ & $100$ & $0.94 \pm 1.28$ & $1.55 \pm 0.9$\\
				& $516 \pm 29.2$ & $100$ & $6.85 \pm 7.34$ & $5.80 \pm 4.3$ & $100$ & $1.20 \pm 0.53$ & $1.24 \pm 0.4$ & $100$ & $1.17 \pm 0.48$ & $1.02 \pm 0.1$\\
				\hline
				
				& $4.0 \pm 2.7$ & $100$ & $0.24 \pm 0.28$ & $1.13 \pm 0.4$ & $100$ & $0.22 \pm 0.24$ & $1.07 \pm 0.3$ & $100$ & $0.22 \pm 0.25$ & $1.1 \pm 0.3$\\
				1024 & $137 \pm 16.6$ & $100$ & $4.17 \pm 7.66$ & $2.58 \pm 2.2$ & $100$ & $2.38 \pm 5.05$ & $1.66 \pm 1.1$ & $100$ & $1.88 \pm 2.79$ & $1.6 \pm 1.0$\\
				& $681 \pm 35.8$ & $100$ & $20.55 \pm 18.7$ & $7.82 \pm 4.9$ & $100$ & $2.99 \pm 1.52$ & $1.42 \pm 0.6$ & $100$ & $2.42 \pm 0.76$ & $1.02 \pm 0.1$\\
				\hline
				
				& $4.6 \pm 3.1$ & $100$ & $0.55 \pm 0.63$ & $1.15 \pm 0.4$ & $100$ & $0.50 \pm 0.55$ & $1.07 \pm 0.3$ & $100$ & $0.44 \pm 0.45$ & $1.03 \pm 0.2$\\
				1536 & $212 \pm 18.8$ & $100$ & $37.67 \pm 66.9$ & $3.93 \pm 3.1$ & $100$ & $10.22 \pm 20.0$ & $1.96 \pm 1.5$ & $100$ & $7.28 \pm 13.0$ & $1.78 \pm 1.3$\\
				& $1017 \pm 39.9$ & $68$ & $106.8 \pm 76.5$ & $12.5 \pm 6.5$ & $100$ & $7.89 \pm 3.55$ & $1.33 \pm 0.6$ & $100$ & $7.04 \pm 2.46$ & $1.03 \pm 0.2$\\
				\hline

				& $3.7 \pm 2.6$ & $100$ & $0.77 \pm 0.81$ & $1.06 \pm 0.2$ & $100$ & $0.88 \pm 0.99$ & $1.06 \pm 0.2$ & $100$ & $0.86 \pm 0.95$ & $1.06 \pm 0.2$\\
				2048 & $284 \pm 22.4$ & $55$ & $132.7 \pm 226.7$ & $4.71 \pm 3.9$ & $100$ & $25.31 \pm 50.5$ & $2.09 \pm 1.6$ & $100$ & $20.69 \pm 51.8$ & $1.97 \pm 1.5$\\
				& $1366 \pm 48.6$ & $23$ & $315.5 \pm 172.8$ & $17.3 \pm 7.3$ & $100$ & $21.40 \pm 16.1$ & $1.58 \pm 1.1$ & $100$ & $14.93 \pm 3.10$ & $1.01 \pm 0.1$\\
				\hline

				\hline
				& $3.7 \pm 2.5$ & $100$ & $2.15 \pm 2.45$ & $1.11 \pm 0.4$ & $100$ & $2.16 \pm 2.64$ & $1.08 \pm 0.3$ & $100$ & $1.83 \pm 2.19$ & $1.06 \pm 0.3$\\
				3072 & $410 \pm 28.8$ & $7$ & $740.5 \pm 764.1$ & $6.00 \pm 4.4$ & $35$ & $188.8 \pm 318.5$ & $3.14 \pm 2.2$ & $74$ & $95.21 \pm 207.2$ & $2.39 \pm 1.7$\\
				& $2036 \pm 46.1$ & $4$ & $1452 \pm 458.3$ & $27.0 \pm 5.9$ & $84$ & $86.01 \pm 72.9$ & $2.07 \pm 1.6$ & $100$ & $45.46 \pm 8.52$ & $1.01 \pm 0.1$\\
				\hline

				& $4.3 \pm 2.8$ & $100$ & $3.48 \pm 3.92$ & $1.10 \pm 0.3$ & $100$ & $3.55 \pm 4.24$ & $1.07 \pm 0.3$ & $100$ & $4.08 \pm 4.84$ & $1.07 \pm 0.3$\\
				4096 & $553 \pm 33.6$ & $-$ & $-$ & $-$ & $23$ & $296.9 \pm 415.6$ & $2.96 \pm 2.0$ & $39$ & $172.2 \pm 262.2$ & $2.33 \pm 1.5$\\
				& $2742 \pm 67.0$ & $-$ & $-$ & $-$ & $44$ & $164.9 \pm 165.6$ & $1.84 \pm 1.6$ & $72$ & $101.3 \pm 22.6$ & $1.03 \pm 0.2$\\
				\hline
			\end{tabular}
		}
	\end{table}

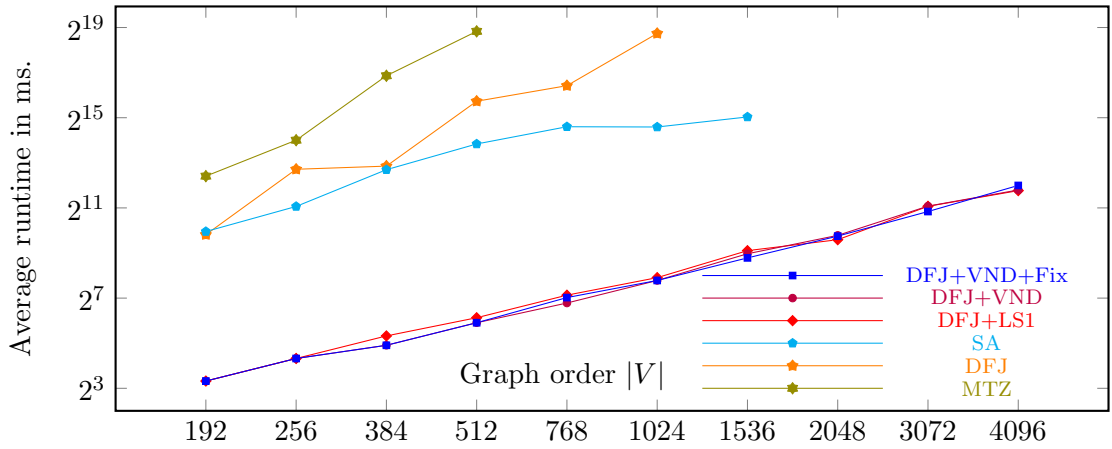
\begin{figure}[p]
	\centering
	\begin{tikzpicture}[scale=1]
	\begin{axis}[
	y=0.3cm,
	x=1.2cm,
	ymode = log,
	log basis y={2},
	axis line style = thick,
	xlabel={Graph order $|V|$},
	x label style={at={(axis description cs:0.45,0.125)},anchor=south},
	ylabel={Average runtime in ms.},
	xtick       = {1,2,3,4,5,6,7,8,9,10},
	xticklabels = {192,256,384,512,768,1024,1536,2048,3072,4096},
	xmin=0,
	xmax=11,
	ymin=4,
	ymax=1000000.]

	\node [diamond,draw,fill,inner sep=1pt,red] (DFJLS1) at (axis cs: 1,10) {};
	\node [diamond,draw,fill,inner sep=1pt,red] (DFJLS2) at (axis cs: 2,20) {};
	\node [diamond,draw,fill,inner sep=1pt,red] (DFJLS3) at (axis cs: 3,40) {};
	\node [diamond,draw,fill,inner sep=1pt,red] (DFJLS4) at (axis cs: 4,70) {};
	\node [diamond,draw,fill,inner sep=1pt,red] (DFJLS5) at (axis cs: 5,140) {};
	\node [diamond,draw,fill,inner sep=1pt,red] (DFJLS6) at (axis cs: 6,240) {};
	\node [diamond,draw,fill,inner sep=1pt,red] (DFJLS7) at (axis cs: 7,550) {};
	\node [diamond,draw,fill,inner sep=1pt,red] (DFJLS8) at (axis cs: 8,770) {};
	\node [diamond,draw,fill,inner sep=1pt,red] (DFJLS9) at (axis cs: 9,2150) {};
	\node [diamond,draw,fill,inner sep=1pt,red] (DFJLS10) at (axis cs: 10,3480) {};
	
	\draw [red] (DFJLS1) -- (DFJLS2) -- (DFJLS3) -- (DFJLS4) -- (DFJLS5) -- (DFJLS6) -- (DFJLS7) -- (DFJLS8) -- (DFJLS9) -- (DFJLS10);

	\node [circle,draw,fill,inner sep=1pt,purple] (DFJVND1) at (axis cs: 1,10) {};
	\node [circle,draw,fill,inner sep=1pt,purple] (DFJVND2) at (axis cs: 2,20) {};
	\node [circle,draw,fill,inner sep=1pt,purple] (DFJVND3) at (axis cs: 3,30) {};
	\node [circle,draw,fill,inner sep=1pt,purple] (DFJVND4) at (axis cs: 4,60) {};
	\node [circle,draw,fill,inner sep=1pt,purple] (DFJVND5) at (axis cs: 5,110) {};
	\node [circle,draw,fill,inner sep=1pt,purple] (DFJVND6) at (axis cs: 6,220) {};
	\node [circle,draw,fill,inner sep=1pt,purple] (DFJVND7) at (axis cs: 7,500) {};
	\node [circle,draw,fill,inner sep=1pt,purple] (DFJVND8) at (axis cs: 8,880) {};
	\node [circle,draw,fill,inner sep=1pt,purple] (DFJVND9) at (axis cs: 9,2160) {};
	\node [circle,draw,fill,inner sep=1pt,purple] (DFJVND10) at (axis cs: 10,3550) {};

	\draw [purple] (DFJVND1) -- (DFJVND2) -- (DFJVND3) -- (DFJVND4) -- (DFJVND5) -- (DFJVND6) -- (DFJVND7) -- (DFJVND8) -- (DFJVND9) -- (DFJVND10);

	\node [draw,fill,inner sep=1.25pt,blue] (DFJVNDFix1) at (axis cs: 1,10) {};	
	\node [draw,fill,inner sep=1.25pt,blue] (DFJVNDFix2) at (axis cs: 2,20) {};
	\node [draw,fill,inner sep=1.25pt,blue] (DFJVNDFix3) at (axis cs: 3,30) {};
	\node [draw,fill,inner sep=1.25pt,blue] (DFJVNDFix4) at (axis cs: 4,60) {};
	\node [draw,fill,inner sep=1.25pt,blue] (DFJVNDFix5) at (axis cs: 5,130) {};
	\node [draw,fill,inner sep=1.25pt,blue] (DFJVNDFix6) at (axis cs: 6,220) {};
	\node [draw,fill,inner sep=1.25pt,blue] (DFJVNDFix7) at (axis cs: 7,440) {};
	\node [draw,fill,inner sep=1.25pt,blue] (DFJVNDFix8) at (axis cs: 8,860) {};
	\node [draw,fill,inner sep=1.25pt,blue] (DFJVNDFix9) at (axis cs: 9,1830) {};
	\node [draw,fill,inner sep=1.25pt,blue] (DFJVNDFix10) at (axis cs: 10,4080) {};
	
	\draw [blue] (DFJVNDFix1) -- (DFJVNDFix2) -- (DFJVNDFix3) -- (DFJVNDFix4) -- (DFJVNDFix5) -- (DFJVNDFix6) -- (DFJVNDFix7) -- (DFJVNDFix8) -- (DFJVNDFix9) -- (DFJVNDFix10);
	
	\node [star,draw,fill,inner sep=1pt,orange] (DFJ1) at (axis cs: 1,900) {};
	\node [star,draw,fill,inner sep=1pt,orange] (DFJ2) at (axis cs: 2,6710) {};
	\node [star,draw,fill,inner sep=1pt,orange] (DFJ3) at (axis cs: 3,7390) {};
	\node [star,draw,fill,inner sep=1pt,orange] (DFJ4) at (axis cs: 4,54300) {};
	\node [star,draw,fill,inner sep=1pt,orange] (DFJ5) at (axis cs: 5,87590) {};
	\node [star,draw,fill,inner sep=1pt,orange] (DFJ6) at (axis cs: 6,435900) {};
	
	\draw [orange] (DFJ1) -- (DFJ2) -- (DFJ3) -- (DFJ4) -- (DFJ5) -- (DFJ6);

	\node [star,star points=6,draw,fill,inner sep=1pt,olive] (MTZ1) at (axis cs: 1,5430) {};
	\node [star,star points=6,draw,fill,inner sep=1pt,olive] (MTZ2) at (axis cs: 2,16400) {};
	\node [star,star points=6,draw,fill,inner sep=1pt,olive] (MTZ3) at (axis cs: 3,118800) {};
	\node [star,star points=6,draw,fill,inner sep=1pt,olive] (MTZ4) at (axis cs: 4,465300) {};
	
	\draw [olive] (MTZ1) -- (MTZ2) -- (MTZ3) -- (MTZ4);

	\node [regular polygon,regular polygon sides=5,draw,fill,inner sep=1pt,cyan] (SA1) at (axis cs: 1,990) {};
	\node [regular polygon,regular polygon sides=5,draw,fill,inner sep=1pt,cyan] (SA2) at (axis cs: 2,2140) {};
	\node [regular polygon,regular polygon sides=5,draw,fill,inner sep=1pt,cyan] (SA3) at (axis cs: 3,6620) {};
	\node [regular polygon,regular polygon sides=5,draw,fill,inner sep=1pt,cyan] (SA4) at (axis cs: 4,14600) {};
	\node [regular polygon,regular polygon sides=5,draw,fill,inner sep=1pt,cyan] (SA5) at (axis cs: 5,24800) {};
	\node [regular polygon,regular polygon sides=5,draw,fill,inner sep=1pt,cyan] (SA6) at (axis cs: 6,24600) {};
	\node [regular polygon,regular polygon sides=5,draw,fill,inner sep=1pt,cyan] (SA7) at (axis cs: 7,33500) {};
	
	\draw [cyan] (SA1) -- (SA2) -- (SA3) -- (SA4) -- (SA5) -- (SA6) -- (SA7);

	\draw [olive] (axis cs: 6.5,8) -- (axis cs: 8.5,8);
	\node [star,star points=6,draw,fill,inner sep=1pt,olive] at (axis cs: 7.5,8) {};
	\node [olive] at (axis cs: 9.65,8) {\footnotesize MTZ};
	
	\draw [orange] (axis cs: 6.5,16) -- (axis cs: 8.5,16);
	\node [star,draw,fill,inner sep=1pt,orange] at (axis cs: 7.5,16) {};
	\node [orange] at (axis cs: 9.65,16) {\footnotesize DFJ};
	
	\draw [cyan] (axis cs: 6.5,32) -- (axis cs: 8.5,32);
	\node [regular polygon,regular polygon sides=5,draw,fill,inner sep=1pt,cyan] at (axis cs: 7.5,32) {};
	\node [cyan] at (axis cs: 9.65,32) {\footnotesize SA};
	
	\draw [red] (axis cs: 6.5,64) -- (axis cs: 8.5,64);
	\node [diamond,draw,fill,inner sep=1pt,red] at (axis cs: 7.5,64) {};
	\node [red] at (axis cs: 9.65,64) {\footnotesize DFJ+LS1};
	
	\draw [purple] (axis cs: 6.5,128) -- (axis cs: 8.5,128);
	\node [circle,draw,fill,inner sep=1pt,purple] at (axis cs: 7.5,128) {};
	\node [purple] at (axis cs: 9.65,128) {\footnotesize DFJ+VND};
	
	\draw [blue] (axis cs: 6.5,256) -- (axis cs: 8.5,256);
	\node [draw,fill,inner sep=1.25pt,blue] at (axis cs: 7.5,256) {};
	\node [blue] at (axis cs: 9.65,256) {\footnotesize DFJ+VND+Fix};
	
	\end{axis}
	\end{tikzpicture}
	\caption{Computational results for random undirected permutations}
	\label{image:results_random_permutation_undirected}
\end{figure}

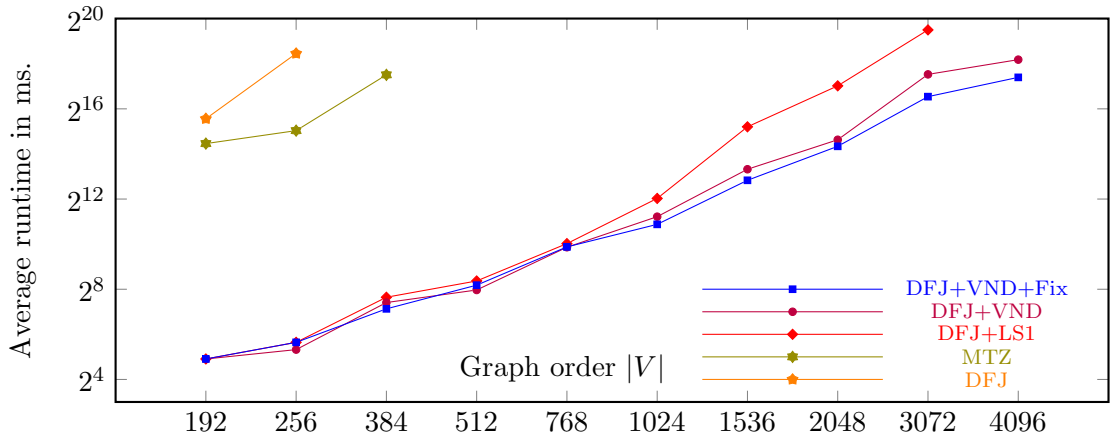
\begin{figure}[p]
	\centering
	\begin{tikzpicture}[scale=1]
	\begin{axis}[
	y=0.3cm,
	x=1.2cm,
	ymode = log,
	log basis y={2},
	axis line style = thick,
	xlabel={Graph order $|V|$},
	x label style={at={(axis description cs:0.45,0.125)},anchor=south},
	ylabel={Average runtime in ms.},
	xtick       = {1,2,3,4,5,6,7,8,9,10},
	xticklabels = {192,256,384,512,768,1024,1536,2048,3072,4096},
	xmin=0,
	xmax=11,
	ymin=8,
	ymax=1050000.]

	\node [diamond,draw,fill,inner sep=1pt,red] (DFJLS1) at (axis cs: 1,30) {};
	\node [diamond,draw,fill,inner sep=1pt,red] (DFJLS2) at (axis cs: 2,50) {};
	\node [diamond,draw,fill,inner sep=1pt,red] (DFJLS3) at (axis cs: 3,200) {};
	\node [diamond,draw,fill,inner sep=1pt,red] (DFJLS4) at (axis cs: 4,330) {};
	\node [diamond,draw,fill,inner sep=1pt,red] (DFJLS5) at (axis cs: 5,1040) {};
	\node [diamond,draw,fill,inner sep=1pt,red] (DFJLS6) at (axis cs: 6,4170) {};
	\node [diamond,draw,fill,inner sep=1pt,red] (DFJLS7) at (axis cs: 7,37670) {};
	\node [diamond,draw,fill,inner sep=1pt,red] (DFJLS8) at (axis cs: 8,132700) {};
	\node [diamond,draw,fill,inner sep=1pt,red] (DFJLS9) at (axis cs: 9,740500) {};
	
	\draw [red] (DFJLS1) -- (DFJLS2) -- (DFJLS3) -- (DFJLS4) -- (DFJLS5) -- (DFJLS6) -- (DFJLS7) -- (DFJLS8)  -- (DFJLS9);

	\node [circle,draw,fill,inner sep=1pt,purple] (DFJVND1) at (axis cs: 1,30) {};
	\node [circle,draw,fill,inner sep=1pt,purple] (DFJVND2) at (axis cs: 2,40) {};
	\node [circle,draw,fill,inner sep=1pt,purple] (DFJVND3) at (axis cs: 3,170) {};
	\node [circle,draw,fill,inner sep=1pt,purple] (DFJVND4) at (axis cs: 4,250) {};
	\node [circle,draw,fill,inner sep=1pt,purple] (DFJVND5) at (axis cs: 5,930) {};
	\node [circle,draw,fill,inner sep=1pt,purple] (DFJVND6) at (axis cs: 6,2380) {};
	\node [circle,draw,fill,inner sep=1pt,purple] (DFJVND7) at (axis cs: 7,10220) {};
	\node [circle,draw,fill,inner sep=1pt,purple] (DFJVND8) at (axis cs: 8,25310) {};
	\node [circle,draw,fill,inner sep=1pt,purple] (DFJVND9) at (axis cs: 9,188800) {};
	\node [circle,draw,fill,inner sep=1pt,purple] (DFJVND10) at (axis cs: 10,296900) {};

	\draw [purple] (DFJVND1) -- (DFJVND2) -- (DFJVND3) -- (DFJVND4) -- (DFJVND5) -- (DFJVND6) -- (DFJVND7) -- (DFJVND8) -- (DFJVND9) -- (DFJVND10);

	\node [draw,fill,inner sep=1.25pt,blue] (DFJVNDFix1) at (axis cs: 1,30) {};	
	\node [draw,fill,inner sep=1.25pt,blue] (DFJVNDFix2) at (axis cs: 2,50) {};
	\node [draw,fill,inner sep=1.25pt,blue] (DFJVNDFix3) at (axis cs: 3,140) {};
	\node [draw,fill,inner sep=1.25pt,blue] (DFJVNDFix4) at (axis cs: 4,290) {};
	\node [draw,fill,inner sep=1.25pt,blue] (DFJVNDFix5) at (axis cs: 5,940) {};
	\node [draw,fill,inner sep=1.25pt,blue] (DFJVNDFix6) at (axis cs: 6,1880) {};
	\node [draw,fill,inner sep=1.25pt,blue] (DFJVNDFix7) at (axis cs: 7,7280) {};
	\node [draw,fill,inner sep=1.25pt,blue] (DFJVNDFix8) at (axis cs: 8,20690) {};
	\node [draw,fill,inner sep=1.25pt,blue] (DFJVNDFix9) at (axis cs: 9,95210) {};
	\node [draw,fill,inner sep=1.25pt,blue] (DFJVNDFix10) at (axis cs: 10,172220) {};
	
	\draw [blue] (DFJVNDFix1) -- (DFJVNDFix2) -- (DFJVNDFix3) -- (DFJVNDFix4) -- (DFJVNDFix5) -- (DFJVNDFix6) -- (DFJVNDFix7) -- (DFJVNDFix8) -- (DFJVNDFix9) -- (DFJVNDFix10);
	
	\node [star,draw,fill,inner sep=1pt,orange] (DFJ1) at (axis cs: 1,48270) {};
	\node [star,draw,fill,inner sep=1pt,orange] (DFJ2) at (axis cs: 2,358600) {};
	
	\draw [orange] (DFJ1) -- (DFJ2);

	\node [star,star points=6,draw,fill,inner sep=1pt,olive] (MTZ1) at (axis cs: 1,22500) {};
	\node [star,star points=6,draw,fill,inner sep=1pt,olive] (MTZ2) at (axis cs: 2,33500) {};
	\node [star,star points=6,draw,fill,inner sep=1pt,olive] (MTZ3) at (axis cs: 3,186100) {};

	\draw [olive] (MTZ1) -- (MTZ2) -- (MTZ3);

	\draw [orange] (axis cs: 6.5,16) -- (axis cs: 8.5,16);
	\node [star,draw,fill,inner sep=1pt,orange] at (axis cs: 7.5,16) {};
	\node [orange] at (axis cs: 9.65,16) {\footnotesize DFJ};
	
	\draw [olive] (axis cs: 6.5,32) -- (axis cs: 8.5,32);
	\node [star,star points=6,draw,fill,inner sep=1pt,olive] at (axis cs: 7.5,32) {};
	\node [olive] at (axis cs: 9.65,32) {\footnotesize MTZ};
	
	\draw [red] (axis cs: 6.5,64) -- (axis cs: 8.5,64);
	\node [diamond,draw,fill,inner sep=1pt,red] at (axis cs: 7.5,64) {};
	\node [red] at (axis cs: 9.65,64) {\footnotesize DFJ+LS1};
	
	\draw [purple] (axis cs: 6.5,128) -- (axis cs: 8.5,128);
	\node [circle,draw,fill,inner sep=1pt,purple] at (axis cs: 7.5,128) {};
	\node [purple] at (axis cs: 9.65,128) {\footnotesize DFJ+VND};
	
	\draw [blue] (axis cs: 6.5,256) -- (axis cs: 8.5,256);
	\node [draw,fill,inner sep=1.25pt,blue] at (axis cs: 7.5,256) {};
	\node [blue] at (axis cs: 9.65,256) {\footnotesize DFJ+VND+Fix};
	
	\end{axis}
	\end{tikzpicture}
	\caption{Computational results for undirected four-peak cycles}
	\label{image:results_undirected_4_peaks}
\end{figure}

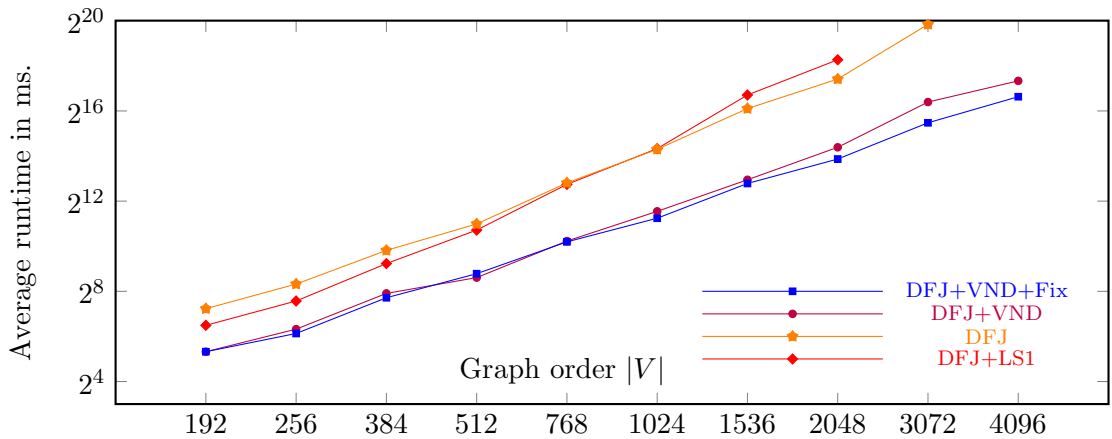
\begin{figure}[p]
	\centering
	\begin{tikzpicture}[scale=1]
	\begin{axis}[
	y=0.3cm,
	x=1.2cm,
	ymode = log,
	log basis y={2},
	axis line style = thick,
	xlabel={Graph order $|V|$},
	x label style={at={(axis description cs:0.45,0.125)},anchor=south},
	ylabel={Average runtime in ms.},
	xtick       = {1,2,3,4,5,6,7,8,9,10},
	xticklabels = {192,256,384,512,768,1024,1536,2048,3072,4096},
	xmin=0,
	xmax=11,
	ymin=8,
	ymax=1050000.]

	\node [diamond,draw,fill,inner sep=1pt,red] (DFJLS1) at (axis cs: 1,90) {};
	\node [diamond,draw,fill,inner sep=1pt,red] (DFJLS2) at (axis cs: 2,190) {};
	\node [diamond,draw,fill,inner sep=1pt,red] (DFJLS3) at (axis cs: 3,600) {};
	\node [diamond,draw,fill,inner sep=1pt,red] (DFJLS4) at (axis cs: 4,1680) {};
	\node [diamond,draw,fill,inner sep=1pt,red] (DFJLS5) at (axis cs: 5,6850) {};
	\node [diamond,draw,fill,inner sep=1pt,red] (DFJLS6) at (axis cs: 6,20550) {};
	\node [diamond,draw,fill,inner sep=1pt,red] (DFJLS7) at (axis cs: 7,106800) {};
	\node [diamond,draw,fill,inner sep=1pt,red] (DFJLS8) at (axis cs: 8,315500) {};
	
	\draw [red] (DFJLS1) -- (DFJLS2) -- (DFJLS3) -- (DFJLS4) -- (DFJLS5) -- (DFJLS6) -- (DFJLS7) -- (DFJLS8);

	\node [circle,draw,fill,inner sep=1pt,purple] (DFJVND1) at (axis cs: 1,40) {};
	\node [circle,draw,fill,inner sep=1pt,purple] (DFJVND2) at (axis cs: 2,80) {};
	\node [circle,draw,fill,inner sep=1pt,purple] (DFJVND3) at (axis cs: 3,240) {};
	\node [circle,draw,fill,inner sep=1pt,purple] (DFJVND4) at (axis cs: 4,390) {};
	\node [circle,draw,fill,inner sep=1pt,purple] (DFJVND5) at (axis cs: 5,1200) {};
	\node [circle,draw,fill,inner sep=1pt,purple] (DFJVND6) at (axis cs: 6,2990) {};
	\node [circle,draw,fill,inner sep=1pt,purple] (DFJVND7) at (axis cs: 7,7890) {};
	\node [circle,draw,fill,inner sep=1pt,purple] (DFJVND8) at (axis cs: 8,21400) {};
	\node [circle,draw,fill,inner sep=1pt,purple] (DFJVND9) at (axis cs: 9,86000) {};
	\node [circle,draw,fill,inner sep=1pt,purple] (DFJVND10) at (axis cs: 10,164900) {};

	\draw [purple] (DFJVND1) -- (DFJVND2) -- (DFJVND3) -- (DFJVND4) -- (DFJVND5) -- (DFJVND6) -- (DFJVND7) -- (DFJVND8) -- (DFJVND9) -- (DFJVND10);

	\node [draw,fill,inner sep=1.25pt,blue] (DFJVNDFix1) at (axis cs: 1,40) {};	
	\node [draw,fill,inner sep=1.25pt,blue] (DFJVNDFix2) at (axis cs: 2,70) {};
	\node [draw,fill,inner sep=1.25pt,blue] (DFJVNDFix3) at (axis cs: 3,210) {};
	\node [draw,fill,inner sep=1.25pt,blue] (DFJVNDFix4) at (axis cs: 4,440) {};
	\node [draw,fill,inner sep=1.25pt,blue] (DFJVNDFix5) at (axis cs: 5,1170) {};
	\node [draw,fill,inner sep=1.25pt,blue] (DFJVNDFix6) at (axis cs: 6,2420) {};
	\node [draw,fill,inner sep=1.25pt,blue] (DFJVNDFix7) at (axis cs: 7,7040) {};
	\node [draw,fill,inner sep=1.25pt,blue] (DFJVNDFix8) at (axis cs: 8,14930) {};
	\node [draw,fill,inner sep=1.25pt,blue] (DFJVNDFix9) at (axis cs: 9,45460) {};
	\node [draw,fill,inner sep=1.25pt,blue] (DFJVNDFix10) at (axis cs: 10,101300) {};
	
	\draw [blue] (DFJVNDFix1) -- (DFJVNDFix2) -- (DFJVNDFix3) -- (DFJVNDFix4) -- (DFJVNDFix5) -- (DFJVNDFix6) -- (DFJVNDFix7) -- (DFJVNDFix8) -- (DFJVNDFix9) -- (DFJVNDFix10);

	\node [star,draw,fill,inner sep=1pt,orange] (DFJ1) at (axis cs: 1,150) {};
	\node [star,draw,fill,inner sep=1pt,orange] (DFJ2) at (axis cs: 2,320) {};
	\node [star,draw,fill,inner sep=1pt,orange] (DFJ3) at (axis cs: 3,900) {};
	\node [star,draw,fill,inner sep=1pt,orange] (DFJ4) at (axis cs: 4,2030) {};
	\node [star,draw,fill,inner sep=1pt,orange] (DFJ5) at (axis cs: 5,7160) {};
	\node [star,draw,fill,inner sep=1pt,orange] (DFJ6) at (axis cs: 6,20130) {};
	\node [star,draw,fill,inner sep=1pt,orange] (DFJ7) at (axis cs: 7,70340) {};
	\node [star,draw,fill,inner sep=1pt,orange] (DFJ8) at (axis cs: 8,174300) {};
	\node [star,draw,fill,inner sep=1pt,orange] (DFJ9) at (axis cs: 9,935100) {};
	
	\draw [orange] (DFJ1) -- (DFJ2) -- (DFJ3) -- (DFJ4) -- (DFJ5) -- (DFJ6) -- (DFJ7) -- (DFJ8) -- (DFJ9);

	\draw [red] (axis cs: 6.5,32) -- (axis cs: 8.5,32);
	\node [diamond,draw,fill,inner sep=1pt,red] at (axis cs: 7.5,32) {};
	\node [red] at (axis cs: 9.65,32) {\footnotesize DFJ+LS1};
	
	\draw [orange] (axis cs: 6.5,64) -- (axis cs: 8.5,64);
	\node [star,draw,fill,inner sep=1pt,orange] at (axis cs: 7.5,64) {};
	\node [orange] at (axis cs: 9.65,64) {\footnotesize DFJ};
	
	\draw [purple] (axis cs: 6.5,128) -- (axis cs: 8.5,128);
	\node [circle,draw,fill,inner sep=1pt,purple] at (axis cs: 7.5,128) {};
	\node [purple] at (axis cs: 9.65,128) {\footnotesize DFJ+VND};
	
	\draw [blue] (axis cs: 6.5,256) -- (axis cs: 8.5,256);
	\node [draw,fill,inner sep=1.25pt,blue] at (axis cs: 7.5,256) {};
	\node [blue] at (axis cs: 9.65,256) {\footnotesize DFJ+VND+Fix};

	\end{axis}
	\end{tikzpicture}
	\caption{Computational results for undirected pyramidal tours}
	\label{image:results_undirected_pyramidal}
\end{figure}

For random permutations (Fig.~\ref{image:results_random_permutation_undirected}, Table~\ref{table:undirected}, line 1) DFJ+LS1, DFJ+VND, and DFJ+VND+Fix showed the best results, solving all $1\,000$ instances and being statistically indistinguishable. It can be concluded that the first neighbourhood structure is sufficient to solve such problems. The rest of the heuristics neither improve nor worsen the algorithm.
Note that random undirected 4-regular multigraphs contain a large number of different Hamiltonian decompositions.
This allows one to find the decomposition by randomized algorithms with polynomial expected running time by random matchings \cite{Kim2001}. In particular, the SA algorithm, which constructs the 2-factors through reduction to random perfect matchings, showed the fourth result, solving $540$ instances out of $1\,000$ and being on average $131.5 \pm 48.8$ times slower than the winners which solved the problem in just 1-2 iterations.
The Dantzig-Fulkerson-Johnson (DFJ) formulation showed the fifth result, solving $492$ instances and lagging on average $2.5 \pm 1.2$ times from SA.
It can be seen that the construction of perfect matchings on undirected graphs is much cheaper than calling the ILP-solver.
Finally, the Miller-Tucker-Zemlin (MTZ) formulation showed the last result, solving only $276$ test problems and lagging  behind DFJ by an average of $8.26 \pm 5.0$ times.

On four-peak cycles (Fig.~\ref{image:results_undirected_4_peaks}, Table~\ref{table:undirected}, line 2) up to $768$ vertices, three algorithms DFJ+LS1, DFJ+VND, DFJ+VND+Fix kept the same level. However, starting from $1\,024$ vertices, their results diverged. DFJ+VND+Fix became the champion, solving $913$ instances out of $1\,000$. The second place is taken by DFJ+VND, which solved $858$ problems and lagged behind DFJ+VND+Fix by an average of $1.5 \pm 0.3$ times on graphs with more than $1\,000$ vertices. DFJ+LS1 took third place, solving $762$ instances and lagging behind DFJ+VND by $3.6 \pm 1.2$ times on large graphs. The fourth place was taken by MTZ, which solved $161$ instances and lost to DFJ+LS1 by $479.8 \pm 324.3$ times. In fifth place, DFJ solved $120$ problems and unexpectedly lost to MTZ by $15.4 \pm 6.5$ times. The SA heuristic showed the last result, solving only $19$ problems.

On pyramidal tours (Fig.~\ref{image:results_undirected_pyramidal}, Table~\ref{table:undirected}, line 3), DFJ+VND+Fix showed the best result, solving $972$ instances. DFJ+VND took second place with $928$ instances and lagged behind DFJ+VND+Fix by an average of $1.25 \pm 0.3$ times. The third place was taken by DFJ with $748$ instances and lagging behind DFJ+VND by an average of $6.37 \pm 2.4$ times. DFJ+LS1 solved $695$ problems and showed results comparable to DFJ, winning on small graphs up to $1\,000$ vertices and losing on large ones. It can be concluded that with a large number of multiple edges, the first neighbourhood structures becomes not very efficient. In the tail are MTZ, which solved only 1 instance, and SA, which did not solve any test instances on pyramidal tours at all.

In general, the DFJ+VND+Fix algorithm showed the best results on undirected graphs. Both neighbourhood structures provided significant performance gains. The chain edge fixing procedure did not slow down the algorithm, and on graphs with a significant number of multiple edges it enhances the performance. The SA algorithm based on random perfect matchings turned out to be more or less efficient only on random permutations and did not work completely on other types of graphs. Of the two ILP-models, the DFJ performed better, except for four-peak cycles, while the MTZ model was too cumbersome and ineffective with a large number of multiple edges.

The results of computational experiments on directed graphs are given in Table~\ref{table:directed}.
The fundamental difference between directed and undirected graphs is the fact that not all instances of the problem had a solution. Among the random permutations, there were only $194$ feasible problems out of $1\,000$, among the four-peak cycles -- $934$, and among the pyramidal tours -- all $1\,000$.

\begin{table}[p]
	\centering
	\caption{Computational results for directed graphs}
	\label{table:directed}
	\resizebox{\textwidth}{!}{
		\begin{tabular}{|*{14}{r|}}
			\hline
			& & \multicolumn{6}{c|}{MTZ} & \multicolumn{6}{c|}{DFJ} \\ 
			\hline
			& & \multicolumn{3}{c|}{Feasible} & \multicolumn{3}{c|}{Infeasible} & \multicolumn{3}{c|}{Feasible} & \multicolumn{3}{c|}{Infeasible} \\ \hline
			$|V|$ & $|x \cap y|$ & N & time (s) & Iter & N & time (s) & Iter & N & time (s) & Iter & N & time (s) & Iter \\ 
			\hline
			& $2.2 \pm 2.0$ & $21$ & $0.03 \pm 0.03$ & $1$ & $79$ & $0.15 \pm 0.11$ & $1$ & $21$ & $0.01 \pm 0.01$ & $3.52 \pm 2.1$ & $79$ & $0.01 \pm 0.01$ & $4.40 \pm 1.9$ \\
			192 & $13.3 \pm 4.5$ & $69$ & $0.09 \pm 0.11$ & $1$ & $31$ & $0.25 \pm 0.10$ & $1$ & $69$ & $0.02 \pm 0.02$ & $7.14 \pm 4.3$ & $31$ & $0.02 \pm 0.01$ & $6.74 \pm 2.6$\\
			& $67.4 \pm 12.8$ & $100$ & $0.17 \pm 0.24$ & $1$ & $-$ & $-$ & $-$ & $100$ & $0.01 \pm 0.00$ & $4.93 \pm 1.5$ & $-$ & $-$ & $-$ \\
			
			\hline
			& $1.9 \pm 1.7$ & $25$ & $0.05 \pm 0.09$ & $1$ & $75$ & $0.26 \pm 0.14$ & $1$ & $25$ & $0.08 \pm 0.15$ & $6.32 \pm 5.9$ & $75$ & $0.06 \pm 0.06$ & $6.01 \pm 2.8$ \\
			256 & $18.5 \pm 6.5$ & $80$ & $0.22 \pm 0.23$ & $1$ & $20$ & $0.48 \pm 0.29$ & $1$ & $80$ & $0.14 \pm 0.19$ & $9.66 \pm 5.5$ & $20$ & $0.05 \pm 0.02$ & $6.9 \pm 1.6$\\
			& $86.9 \pm 14.7$ & $100$ & $0.47 \pm 0.52$ & $1$ & $-$ & $-$ & $-$ & $100$ & $0.07 \pm 0.04$ & $5.29 \pm 2.9$ & $-$ & $-$ & $-$\\
			\hline

			& $2.1 \pm 2.1$ & $20$ & $0.10 \pm 0.17$ & $1$ & $80$ & $0.55 \pm 0.36$ & $1$ & $20$ & $0.05 \pm 0.03$ & $4.9 \pm 2.3$ & $80$ & $0.07 \pm 0.11$ & $5.76 \pm 3.8$ \\
			384 & $26.7 \pm 7.4$ & $90$ & $0.76 \pm 0.74$ & $1$ & $10$ & $1.44 \pm 0.88$ & $1$ & $90$ & $0.25 \pm 0.25$ & $11.4 \pm 5.7$ & $10$ & $0.10 \pm 0.06$ & $7.7 \pm 1.6$\\
			& $129 \pm 20.5$ & $100$ & $1.80 \pm 2.52$ & $1$ & $-$ & $-$ & $-$ & $100$ & $0.07 \pm 0.05$ & $5.76 \pm 3.0$ & $-$ & $-$ & $-$\\
			\hline

			& $1.9 \pm 1.9$ & $22$ & $0.23 \pm 0.38$ & $1$ & $78$ & $0.91 \pm 0.64$ & $1$ & $22$ & $0.05 \pm 0.03$ & $5.05 \pm 2.6$ & $78$ & $0.07 \pm 0.13$ & $5.88 \pm 4.0$ \\
			512 & $35.0 \pm 8.9$ & $96$ & $1.60 \pm 1.45$ & $1$ & $4$ & $5.16 \pm 1.39$ & $1$ & $96$ & $0.27 \pm 0.29$ & $13.3 \pm 7.4$ & $4$ & $0.20 \pm 0.16$ & $12.0 \pm 4.6$\\
			& $171 \pm 23.5$ & $100$ & $8.45 \pm 11.8$ & $1$ & $-$ & $-$ & $-$ & $100$ & $0.09 \pm 0.05$ & $6.42 \pm 3.3$ & $-$ & $-$ & $-$\\
			\hline

			& $2.0 \pm 2.0$ & $19$ & $1.08 \pm 1.00$ & $1$ & $81$ & $2.14 \pm 1.07$ & $1$ & $19$ & $0.11 \pm 0.09$ & $6.21 \pm 4.1$ & $81$ & $0.09 \pm 0.06$ & $5.67 \pm 2.3$ \\
			768 & $53.7 \pm 10.8$ & $99$ & $4.84 \pm 3.44$ & $1$ & $1$ & $13.2 \pm 0.0$ & $1$ & $99$ & $1.45 \pm 2.12$ & $23.4 \pm 15.0$ & $1$ & $0.32 \pm 0.00$ & $12 \pm 0.0$\\
			& $258 \pm 26.3$ & $100$ & $60.9 \pm 88.8$ & $1$ & $-$ & $-$ & $-$ & $100$ & $0.15 \pm 0.11$ & $7.08 \pm 4.6$ & $-$ & $-$ & $-$\\
			\hline

			& $2.0 \pm 1.9$ & $17$ & $2.17 \pm 1.14$ & $1$ & $83$ & $2.73 \pm 1.34$ & $1$ & $17$ & $0.18 \pm 0.17$ & $5.76 \pm 3.6$ & $83$ & $0.20 \pm 0.23$ & $6.46 \pm 3.5$ \\
			1024 & $66.5 \pm 12.3$ & $100$ & $8.57 \pm 5.66$ & $1$ & $-$ & $-$ & $-$ & $100$ & $2.94 \pm 3.68$ & $33.7 \pm 22.1$ & $-$ & $-$ & $-$\\
			& $341 \pm 33.7$ & $28$ & $165.2 \pm 155$ & $1$ & $-$ & $-$ & $-$ & $100$ & $0.18 \pm 0.15$ & $8.13 \pm 5.9$ & $-$ & $-$ & $-$\\
			\hline

			& $2.2 \pm 2.0$ & $16$ & $5.47 \pm 2.01$ & $1$ & $84$ & $4.88 \pm 2.15$ & $1$ & $16$ & $0.14 \pm 0.09$ & $6.5 \pm 3.1$ & $84$ & $0.17 \pm 0.36$ & $7.02 \pm 4.1$ \\
			1536 & $102 \pm 14.0$ & $100$ & $35.5 \pm 61.8$ & $1$ & $-$ & $-$ & $-$ & $100$ & $17.35 \pm 25.7$ & $70.3 \pm 45.1$ & $-$ & $-$ & $-$\\
			& $508 \pm 44.3$ & $0$ & $-$ & $-$ & $-$ & $-$ & $-$ & $100$ & $0.45 \pm 0.35$ & $11.6 \pm 7.5$ & $-$ & $-$ & $-$\\
			\hline

			& $1.5 \pm 1.6$ & $15$ & $8.20 \pm 5.13$ & $1$ & $85$ & $8.34 \pm 3.97$ & $1$ & $15$ & $0.28 \pm 0.25$ & $7.33 \pm 4.4$ & $85$ & $0.22 \pm 0.16$ & $7.05 \pm 3.0$ \\
			2048 & $143 \pm 14.8$ & $40$ & $138.7 \pm 258$ & $1$ & $-$ & $-$ & $-$ & $88$ & $80.1 \pm 130$ & $121.9 \pm 88.2$ & $-$ & $-$ & $-$\\
			& $678 \pm 49.6$ & $-$ & $-$ & $-$ & $-$ & $-$ & $-$ & $100$ & $0.91 \pm 0.73$ & $13.4 \pm 9.4$ & $-$ & $-$ & $-$\\
			\hline

			& $1.9 \pm 2.0$ & $21$ & $18.2 \pm 9.63$ & $1$ & $79$ & $17.8 \pm 9.50$ & $-$ & $21$ & $0.73 \pm 1.52$ & $8.00 \pm 7.5$ & $79$ & $0.48 \pm 0.46$ & $7.91 \pm 3.9$ \\
			3072 & $200 \pm 19.0$ & $24$ & $309.3 \pm 339$ & $-$ & $-$ & $-$ & $-$ & $10$ & $566.7 \pm 371$ & $283.6 \pm 108$ & $-$ & $-$ & $-$\\
			& $1023 \pm 48.9$ & $-$ & $-$ & $-$ & $-$ & $-$ & $-$ & $100$ & $2.80 \pm 2.28$ & $19.5 \pm 13.6$ & $-$ & $-$ & $-$\\
			\hline

			& $2.2 \pm 2.0$ & $18$ & $30.7 \pm 23.4$ & $1$ & $82$ & $32.1 \pm 17.4$ & $1$ & $18$ & $0.77 \pm 0.94$ & $7.89 \pm 4.8$ & $82$ & $0.86 \pm 1.22$ & $8.26 \pm 4.7$ \\
			4096 & $278 \pm 22.8$ & $-$ & $-$ & $-$ & $-$ & $-$ & $-$ & $-$ & $-$ & $-$ & $-$ & $-$ & $-$\\
			& $1377 \pm 71.1$ & $-$ & $-$ & $-$ & $-$ & $-$ & $-$ & $100$ & $5.87 \pm 6.29$ & $20.1 \pm 19.2$ & $-$ & $-$ & $-$\\
			\hline
			\hline

			& & \multicolumn{6}{c|}{SA} & \multicolumn{6}{c|}{DFJ + LS} \\ 
			\hline
			& & \multicolumn{3}{c|}{Solved} & \multicolumn{3}{c|}{Not solved} & \multicolumn{3}{c|}{Feasible} & \multicolumn{3}{c|}{Infeasible} \\ \hline
			$|V|$ & $|x \cap y|$ & N & time (s) & Iter & N & time (s) & Iter & N & time (s) & Iter & N & time (s) & Iter \\ 
			\hline
			& $2.2 \pm 2.0$ & $19$ & $0.54 \pm 0.73$ & $144 \pm 202$ & $81$ & $8.89 \pm 1.34$ & $2500$ & $21$ & $0.03 \pm 0.03$ & $2.14 \pm 1.6$ & $79$ & $0.05 \pm 0.02$ & $3.57 \pm 1.5$ \\
			192 & $13.3 \pm 4.5$ & $42$ & $1.57 \pm 1.87$ & $543 \pm 659$ & $58$ & $7.15 \pm 0.86$ & $2500$ & $69$ & $0.04 \pm 0.03$ & $3.3 \pm 2.2$ & $31$ & $0.07 \pm 0.03$ & $5.42 \pm 1.9$\\
			& $67.4 \pm 12.8$ & $58$ & $0.41 \pm 1.20$ & $149 \pm 499$ & $42$ & $5.99 \pm 0.68$ & $2500$ & $100$ & $0.01 \pm 0.01$ & $1.28 \pm 0.7$ & $-$ & $-$ & $-$\\
			
			\hline
			& $1.9 \pm 1.7$ & $19$ & $1.60 \pm 2.35$ & $316 \pm 477$ & $81$ & $12.4 \pm 1.16$ & $2500$ & $25$ & $0.07 \pm 0.10$ & $3.68 \pm 4.1$ & $75$ & $0.09 \pm 0.05$ & $4.61 \pm 2.1$ \\
			256 & $18.5 \pm 6.5$ & $37$ & $4.37 \pm 2.98$ & $945 \pm 654$ & $63$ & $11.0 \pm 0.45$ & $2500$ & $80$ & $0.09 \pm 0.08$ & $4.6 \pm 3.2$ & $20$ & $0.11 \pm 0.04$ & $5.25 \pm 1.5$\\
			& $86.9 \pm 14.7$ & $37$ & $0.08 \pm 0.04$ & $18.5 \pm 9.8$ & $63$ & $9.18 \pm 0.33$ & $2500$ & $100$ & $0.02 \pm 0.02$ & $1.44 \pm 0.7$ & $-$ & $-$ & $-$\\
			\hline

			& $2.1 \pm 2.1$ & $17$ & $3.65 \pm 4.36$ & $326 \pm 392$ & $84$ & $22.5 \pm 9.10$ & $2500$ & $20$ & $0.09 \pm 0.08$ & $2.4 \pm 1.3$ & $80$ & $0.20 \pm 0.16$ & $4.55 \pm 2.8$ \\
			384 & $26.7 \pm 7.4$ & $16$ & $13.6 \pm 7.40$ & $1341 \pm 778$ & $84$ & $23.7 \pm 3.78$ & $2500$ & $90$ & $0.19 \pm 0.15$ & $4.93 \pm 2.8$ & $10$ & $0.25 \pm 0.12$ & $6.1 \pm 1.4$\\
			& $129 \pm 20.5$ & $18$ & $0.19 \pm 0.06$ & $23.9 \pm 7.9$ & $82$ & $19.2 \pm 3.24$ & $2500$ & $100$ & $0.04 \pm 0.05$ & $1.6 \pm 1.1$ & $-$ & $-$ & $-$\\
			\hline

			& $1.9 \pm 1.9$ & $17$ & $11.7 \pm 9.64$ & $576 \pm 501$ & $83$ & $46.5 \pm 5.36$ & $2500$ & $22$ & $0.24 \pm 0.19$ & $3.04 \pm 2.0$ & $78$ & $0.36 \pm 0.25$ & $4.67 \pm 3.1$ \\
			512 & $35.0 \pm 8.9$ & $-$ & $-$ & $-$ & $-$ & $-$ & $-$ & $96$ & $0.34 \pm 0.27$ & $5.20 \pm 3.4$ & $4$ & $0.62 \pm 0.24$ & $9.5 \pm 2.1$\\
			& $171 \pm 23.5$ & $-$ & $-$ & $-$ & $-$ & $-$ & $-$ & $100$ & $0.06 \pm 0.08$ & $1.5 \pm 0.9$ & $-$ & $-$ & $-$\\
			\hline

			& $2.0 \pm 2.0$ & $11$ & $16.5 \pm 16.8$ & $363 \pm 438$ & $89$ & $60.0$ & $1388$ & $19$ & $0.62 \pm 0.50$ & $3.37 \pm 2.5$ & $81$ & $0.81 \pm 0.42$ & $4.44 \pm 1.9$ \\
			768 & $53.7 \pm 10.8$ & $-$ & $-$ & $-$ & $-$ & $-$ & $-$ & $99$ & $0.81 \pm 0.67$ & $6.68 \pm 4.8$ & $1$ & $0.99 \pm 0.00$ & $8 \pm 0.0$\\
			& $258 \pm 26.3$ & $-$ & $-$ & $-$ & $-$ & $-$ & $-$ & $100$ & $0.13 \pm 0.22$ & $1.53 \pm 1.1$ & $-$ & $-$ & $-$\\
			\hline
			
			& $2.0 \pm 1.9$ & $6$ & $15.8 \pm 14.1$ & $177 \pm 167$ & $94$ & $60.0$ & $764$ & $17$ & $0.72 \pm 0.83$ & $2.65 \pm 2.0$ & $83$ & $1.59 \pm 0.88$ & $5.0 \pm 2.6$ \\
			1024 & $66.5 \pm 12.3$ & $-$ & $-$ & $-$ & $-$ & $-$ & $-$ & $100$ & $1.85 \pm 1.85$ & $8.49 \pm 6.8$ & $-$ & $-$ & $-$\\
			& $341.4 \pm 33.7$ & $-$ & $-$ & $-$ & $-$ & $-$ & $-$ & $100$ & $0.27 \pm 0.49$ & $1.69 \pm 1.4$ & $-$ & $-$ & $-$\\
			\hline

			& $2.2 \pm 2.0$ & $2$ & $17.5 \pm 7.94$ & $380 \pm 49.5$ & $98$ & $60.0$ & $381$ & $16$ & $2.43 \pm 1.67$ & $3.25 \pm 1.8$ & $84$ & $4.61 \pm 2.44$ & $5.55 \pm 3.1$ \\
			1536 & $102 \pm 14.0$ & $-$ & $-$ & $-$ & $-$ & $-$ & $-$ & $100$ & $4.65 \pm 5.21$ & $9.33 \pm 8.9$ & $-$ & $-$ & $-$\\
			& $508 \pm 44.3$ & $-$ & $-$ & $-$ & $-$ & $-$ & $-$ & $100$ & $0.49 \pm 1.16$ & $1.51 \pm 1.4$ & $-$ & $-$ & $-$\\
			\hline

			& $1.5 \pm 1.6$ & $-$ & $-$ & $-$ & $-$ & $-$ & $-$ & $15$ & $6.45 \pm 5.32$ & $4.0 \pm 2.6$ & $85$ & $8.35 \pm 3.50$ & $5.41 \pm 2.3$ \\
			2048 & $143 \pm 14.8$ & $-$ & $-$ & $-$ & $-$ & $-$ & $-$ & $100$ & $16.13 \pm 19.1$ & $15.9 \pm 15.2$ & $-$ & $-$ & $-$\\
			& $678 \pm 49.6$ & $-$ & $-$ & $-$ & $-$ & $-$ & $-$ & $100$ & $0.95 \pm 2.20$ & $1.54 \pm 1.4$ & $-$ & $-$ & $-$\\
			\hline

			& $1.9 \pm 2.0$ & $-$ & $-$ & $-$ & $-$ & $-$ & $-$ & $21$ & $8.87 \pm 11.5$ & $2.76 \pm 2.6$ & $79$ & $22.6 \pm 12.1$ & $5.95 \pm 2.7$ \\
			3072 & $200 \pm 19.0$ & $-$ & $-$ & $-$ & $-$ & $-$ & $-$ & $69$ & $103.7 \pm 132$ & $35.7 \pm 36.7$ & $-$ & $-$ & $-$\\
			& $1023 \pm 48.9$ & $-$ & $-$ & $-$ & $-$ & $-$ & $-$ & $100$ & $5.40 \pm 9.14$ & $2.3 \pm 2.3$ & $-$ & $-$ & $-$\\
			\hline

			& $2.2 \pm 2.0$ & $-$ & $-$ & $-$ & $-$ & $-$ & $-$ & $18$ & $32.3 \pm 25.6$ & $4.94 \pm 4.2$ & $82$ & $47.0 \pm 35.5$ & $6.11 \pm 3.5$ \\
			4096 & $278 \pm 22.8$ & $-$ & $-$ & $-$ & $-$ & $-$ & $-$ & $28$ & $248.4 \pm 291$ & $43.5 \pm 44.3$ & $-$ & $-$ & $-$\\
			& $1377 \pm 71.1$ & $-$ & $-$ & $-$ & $-$ & $-$ & $-$ & $100$ & $7.54 \pm 17.8$ & $1.94 \pm 2.3$ & $-$ & $-$ & $-$\\
			\hline
		\end{tabular}
	}
\end{table}

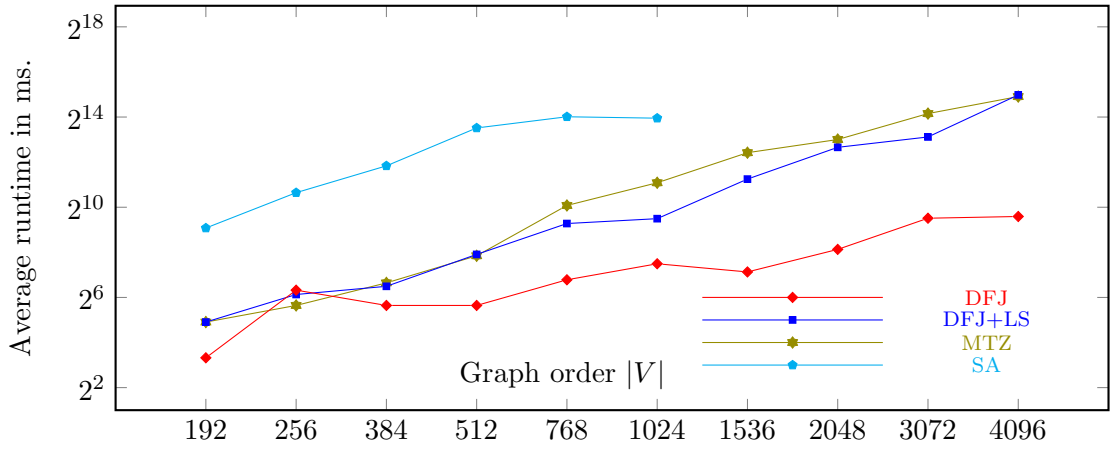
\begin{figure}[p]
	\centering
	\begin{tikzpicture}[scale=1]
	\begin{axis}[
	y=0.3cm,
	x=1.2cm,
	ymode = log,
	log basis y={2},
	axis line style = thick,
	xlabel={Graph order $|V|$},
	x label style={at={(axis description cs:0.45,0.125)},anchor=south},
	ylabel={Average runtime in ms.},
	xtick       = {1,2,3,4,5,6,7,8,9,10},
	xticklabels = {192,256,384,512,768,1024,1536,2048,3072,4096},
	xmin=0,
	xmax=11,
	ymin=2,
	ymax=500000.]

	\node [regular polygon,regular polygon sides=5,draw,fill,inner sep=1pt,cyan] (SA1) at (axis cs: 1,540) {};
	\node [regular polygon,regular polygon sides=5,draw,fill,inner sep=1pt,cyan] (SA2) at (axis cs: 2,1600) {};
	\node [regular polygon,regular polygon sides=5,draw,fill,inner sep=1pt,cyan] (SA3) at (axis cs: 3,3650) {};
	\node [regular polygon,regular polygon sides=5,draw,fill,inner sep=1pt,cyan] (SA4) at (axis cs: 4,11700) {};
	\node [regular polygon,regular polygon sides=5,draw,fill,inner sep=1pt,cyan] (SA5) at (axis cs: 5,16500) {};
	\node [regular polygon,regular polygon sides=5,draw,fill,inner sep=1pt,cyan] (SA6) at (axis cs: 6,15800) {};

	\draw [cyan] (SA1) -- (SA2) -- (SA3) -- (SA4) -- (SA5) -- (SA6);

	\node [star,star points=6,draw,fill,inner sep=1pt,olive] (MTZ1) at (axis cs: 1,30) {};
	\node [star,star points=6,draw,fill,inner sep=1pt,olive] (MTZ2) at (axis cs: 2,50) {};
	\node [star,star points=6,draw,fill,inner sep=1pt,olive] (MTZ3) at (axis cs: 3,100) {};
	\node [star,star points=6,draw,fill,inner sep=1pt,olive] (MTZ4) at (axis cs: 4,230) {};
	\node [star,star points=6,draw,fill,inner sep=1pt,olive] (MTZ5) at (axis cs: 5,1080) {};
	\node [star,star points=6,draw,fill,inner sep=1pt,olive] (MTZ6) at (axis cs: 6,2170) {};
	\node [star,star points=6,draw,fill,inner sep=1pt,olive] (MTZ7) at (axis cs: 7,5470) {};
	\node [star,star points=6,draw,fill,inner sep=1pt,olive] (MTZ8) at (axis cs: 8,8200) {};
	\node [star,star points=6,draw,fill,inner sep=1pt,olive] (MTZ9) at (axis cs: 9,18200) {};
	\node [star,star points=6,draw,fill,inner sep=1pt,olive] (MTZ10) at (axis cs: 10,30700) {};
	
	\draw [olive] (MTZ1) -- (MTZ2) -- (MTZ3) -- (MTZ4) -- (MTZ5) -- (MTZ6) -- (MTZ7) -- (MTZ8) -- (MTZ9) -- (MTZ10);
	
	\node [draw,fill,inner sep=1.25pt,blue] (DFJLS1) at (axis cs: 1,30) {};	
	\node [draw,fill,inner sep=1.25pt,blue] (DFJLS2) at (axis cs: 2,70) {};
	\node [draw,fill,inner sep=1.25pt,blue] (DFJLS3) at (axis cs: 3,90) {};
	\node [draw,fill,inner sep=1.25pt,blue] (DFJLS4) at (axis cs: 4,240) {};
	\node [draw,fill,inner sep=1.25pt,blue] (DFJLS5) at (axis cs: 5,620) {};
	\node [draw,fill,inner sep=1.25pt,blue] (DFJLS6) at (axis cs: 6,720) {};
	\node [draw,fill,inner sep=1.25pt,blue] (DFJLS7) at (axis cs: 7,2430) {};
	\node [draw,fill,inner sep=1.25pt,blue] (DFJLS8) at (axis cs: 8,6450) {};
	\node [draw,fill,inner sep=1.25pt,blue] (DFJLS9) at (axis cs: 9,8870) {};
	\node [draw,fill,inner sep=1.25pt,blue] (DFJLS10) at (axis cs: 10,32300) {};
	
	\draw [blue] (DFJLS1) -- (DFJLS2) -- (DFJLS3) -- (DFJLS4) -- (DFJLS5) -- (DFJLS6) -- (DFJLS7) -- (DFJLS8) -- (DFJLS9) -- (DFJLS10);
	
	\node [diamond,draw,fill,inner sep=1pt,red] (DFJ1) at (axis cs: 1,10) {};
	\node [diamond,draw,fill,inner sep=1pt,red] (DFJ2) at (axis cs: 2,80) {};
	\node [diamond,draw,fill,inner sep=1pt,red] (DFJ3) at (axis cs: 3,50) {};
	\node [diamond,draw,fill,inner sep=1pt,red] (DFJ4) at (axis cs: 4,50) {};
	\node [diamond,draw,fill,inner sep=1pt,red] (DFJ5) at (axis cs: 5,110) {};
	\node [diamond,draw,fill,inner sep=1pt,red] (DFJ6) at (axis cs: 6,180) {};
	\node [diamond,draw,fill,inner sep=1pt,red] (DFJ7) at (axis cs: 7,140) {};
	\node [diamond,draw,fill,inner sep=1pt,red] (DFJ8) at (axis cs: 8,280) {};
	\node [diamond,draw,fill,inner sep=1pt,red] (DFJ9) at (axis cs: 9,730) {};
	\node [diamond,draw,fill,inner sep=1pt,red] (DFJ10) at (axis cs: 10,770) {};
	
	\draw [red] (DFJ1) -- (DFJ2) -- (DFJ3) -- (DFJ4) -- (DFJ5) -- (DFJ6) -- (DFJ7) -- (DFJ8) -- (DFJ9) -- (DFJ10);

	\draw [red] (axis cs: 6.5,64) -- (axis cs: 8.5,64);
	\node [diamond,draw,fill,inner sep=1pt,red] at (axis cs: 7.5,64) {};
	\node [red] at (axis cs: 9.65,64) {\footnotesize DFJ};

	\draw [blue] (axis cs: 6.5,32) -- (axis cs: 8.5,32);
	\node [draw,fill,inner sep=1.25pt,blue] at (axis cs: 7.5,32) {};
	\node [blue] at (axis cs: 9.65,32) {\footnotesize DFJ+LS};
	
	\draw [olive] (axis cs: 6.5,16) -- (axis cs: 8.5,16);
	\node [star,star points=6,draw,fill,inner sep=1pt,olive] at (axis cs: 7.5,16) {};
	\node [olive] at (axis cs: 9.65,16) {\footnotesize MTZ};
	
	\draw [cyan] (axis cs: 6.5,8) -- (axis cs: 8.5,8);
	\node [regular polygon,regular polygon sides=5,draw,fill,inner sep=1pt,cyan] at (axis cs: 7.5,8) {};
	\node [cyan] at (axis cs: 9.65,8) {\footnotesize SA};
	
	\end{axis}
	\end{tikzpicture}
	\caption{Computational results for random directed permutations (feasible problems)}
	\label{image:results_directed_permutations}
\end{figure}

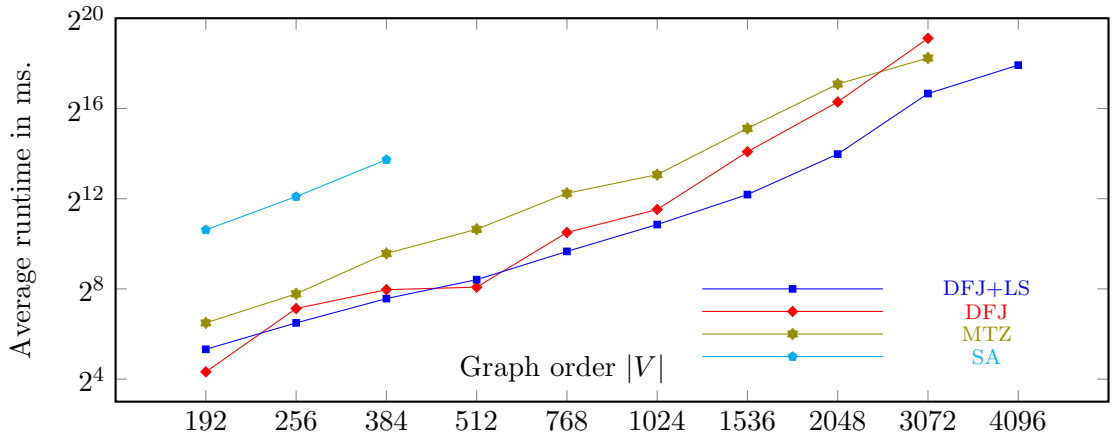
\begin{figure}[p]
	\centering
	\begin{tikzpicture}[scale=1]
	\begin{axis}[
	y=0.3cm,
	x=1.2cm,
	ymode = log,
	log basis y={2},
	axis line style = thick,
	ylabel={Average runtime in ms.},
	xlabel={Graph order $|V|$},
	x label style={at={(axis description cs:0.45,0.125)},anchor=south},
	xtick       = {1,2,3,4,5,6,7,8,9,10},
	xticklabels = {192,256,384,512,768,1024,1536,2048,3072,4096},
	xmin=0,
	xmax=11,
	ymin=8,
	ymax=1050000.]

	\node [star,star points=6,draw,fill,inner sep=1pt,olive] (MTZ1) at (axis cs: 1,90) {};
	\node [star,star points=6,draw,fill,inner sep=1pt,olive] (MTZ2) at (axis cs: 2,220) {};
	\node [star,star points=6,draw,fill,inner sep=1pt,olive] (MTZ3) at (axis cs: 3,760) {};
	\node [star,star points=6,draw,fill,inner sep=1pt,olive] (MTZ4) at (axis cs: 4,1600) {};
	\node [star,star points=6,draw,fill,inner sep=1pt,olive] (MTZ5) at (axis cs: 5,4840) {};
	\node [star,star points=6,draw,fill,inner sep=1pt,olive] (MTZ6) at (axis cs: 6,8570) {};
	\node [star,star points=6,draw,fill,inner sep=1pt,olive] (MTZ7) at (axis cs: 7,35500) {};
	\node [star,star points=6,draw,fill,inner sep=1pt,olive] (MTZ8) at (axis cs: 8,138700) {};
	\node [star,star points=6,draw,fill,inner sep=1pt,olive] (MTZ9) at (axis cs: 9,309300) {};
	
	\draw [olive] (MTZ1) --  (MTZ2) -- (MTZ3) -- (MTZ4) -- (MTZ5) -- (MTZ6) -- (MTZ7) -- (MTZ8) -- (MTZ9);

	\node [regular polygon,regular polygon sides=5,draw,fill,inner sep=1pt,cyan] (SA1) at (axis cs: 1,1570) {};
	\node [regular polygon,regular polygon sides=5,draw,fill,inner sep=1pt,cyan] (SA2) at (axis cs: 2,4370) {};
	\node [regular polygon,regular polygon sides=5,draw,fill,inner sep=1pt,cyan] (SA3) at (axis cs: 3,13600) {};

	\draw [cyan] (SA1) -- (SA2) -- (SA3);
	
	\node [diamond,draw,fill,inner sep=1pt,red] (DFJ1) at (axis cs: 1,20) {};
	\node [diamond,draw,fill,inner sep=1pt,red] (DFJ2) at (axis cs: 2,140) {};
	\node [diamond,draw,fill,inner sep=1pt,red] (DFJ3) at (axis cs: 3,250) {};
	\node [diamond,draw,fill,inner sep=1pt,red] (DFJ4) at (axis cs: 4,270) {};
	\node [diamond,draw,fill,inner sep=1pt,red] (DFJ5) at (axis cs: 5,1450) {};
	\node [diamond,draw,fill,inner sep=1pt,red] (DFJ6) at (axis cs: 6,2940) {};
	\node [diamond,draw,fill,inner sep=1pt,red] (DFJ7) at (axis cs: 7,17350) {};
	\node [diamond,draw,fill,inner sep=1pt,red] (DFJ8) at (axis cs: 8,80100) {};
	\node [diamond,draw,fill,inner sep=1pt,red] (DFJ9) at (axis cs: 9,566700) {};
	
	\draw [red] (DFJ1) -- (DFJ2) -- (DFJ3) -- (DFJ4) -- (DFJ5) -- (DFJ6) -- (DFJ7) -- (DFJ8) -- (DFJ9);

	\node [draw,fill,inner sep=1.25pt,blue] (DFJLS1) at (axis cs: 1,40) {};	
	\node [draw,fill,inner sep=1.25pt,blue] (DFJLS2) at (axis cs: 2,90) {};
	\node [draw,fill,inner sep=1.25pt,blue] (DFJLS3) at (axis cs: 3,190) {};
	\node [draw,fill,inner sep=1.25pt,blue] (DFJLS4) at (axis cs: 4,340) {};
	\node [draw,fill,inner sep=1.25pt,blue] (DFJLS5) at (axis cs: 5,810) {};
	\node [draw,fill,inner sep=1.25pt,blue] (DFJLS6) at (axis cs: 6,1850) {};
	\node [draw,fill,inner sep=1.25pt,blue] (DFJLS7) at (axis cs: 7,4650) {};
	\node [draw,fill,inner sep=1.25pt,blue] (DFJLS8) at (axis cs: 8,16130) {};
	\node [draw,fill,inner sep=1.25pt,blue] (DFJLS9) at (axis cs: 9,103700) {};
	\node [draw,fill,inner sep=1.25pt,blue] (DFJLS10) at (axis cs: 10,248400) {};
	
	\draw [blue] (DFJLS1) -- (DFJLS2) -- (DFJLS3) -- (DFJLS4) -- (DFJLS5) -- (DFJLS6) -- (DFJLS7) -- (DFJLS8) -- (DFJLS9) -- (DFJLS10);

	\draw [blue] (axis cs: 6.5,256) -- (axis cs: 8.5,256);
	\node [draw,fill,inner sep=1.25pt,blue] at (axis cs: 7.5,256) {};
	\node [blue] at (axis cs: 9.65,256) {\footnotesize DFJ+LS};
	
	\draw [red] (axis cs: 6.5,128) -- (axis cs: 8.5,128);
	\node [diamond,draw,fill,inner sep=1pt,red] at (axis cs: 7.5,128) {};
	\node [red] at (axis cs: 9.65,128) {\footnotesize DFJ};
	
	\draw [olive] (axis cs: 6.5,64) -- (axis cs: 8.5,64);
	\node [star,star points=6,draw,fill,inner sep=1pt,olive] at (axis cs: 7.5,64) {};
	\node [olive] at (axis cs: 9.65,64) {\footnotesize MTZ};
	
	\draw [cyan] (axis cs: 6.5,32) -- (axis cs: 8.5,32);
	\node [regular polygon,regular polygon sides=5,draw,fill,inner sep=1pt,cyan] at (axis cs: 7.5,32) {};
	\node [cyan] at (axis cs: 9.65,32) {\footnotesize SA};

	\end{axis}
	\end{tikzpicture}
	\caption{Computational results for directed four-peak cycles (feasible problems)}
	\label{image:results_directed_4_peaks}
\end{figure}

\begin{figure}[p]
	\centering
	\begin{tikzpicture}[scale=1]
	\begin{axis}[
	y=0.3cm,
	x=1.2cm,
	ymode = log,
	log basis y={2},
	axis line style = thick,
	xlabel={Graph order $|V|$},
	x label style={at={(axis description cs:0.45,0.125)},anchor=south},
	ylabel={Average runtime in ms.},
	xtick       = {1,2,3,4,5,6,7,8,9,10},
	xticklabels = {192,256,384,512,768,1024,1536,2048,3072,4096},
	xmin=0,
	xmax=11,
	ymin=2,
	ymax=500000.]

	\node [star,star points=6,draw,fill,inner sep=1pt,olive] (MTZ1) at (axis cs: 1,170) {};
	\node [star,star points=6,draw,fill,inner sep=1pt,olive] (MTZ2) at (axis cs: 2,470) {};
	\node [star,star points=6,draw,fill,inner sep=1pt,olive] (MTZ3) at (axis cs: 3,1800) {};
	\node [star,star points=6,draw,fill,inner sep=1pt,olive] (MTZ4) at (axis cs: 4,8450) {};
	\node [star,star points=6,draw,fill,inner sep=1pt,olive] (MTZ5) at (axis cs: 5,60900) {};
	\node [star,star points=6,draw,fill,inner sep=1pt,olive] (MTZ6) at (axis cs: 6,165200) {};
	
	\draw [olive] (MTZ1) -- (MTZ2) -- (MTZ3) -- (MTZ4) -- (MTZ5) -- (MTZ6);

	\node [regular polygon,regular polygon sides=5,draw,fill,inner sep=1pt,cyan] (SA1) at (axis cs: 1,410) {};
	\node [regular polygon,regular polygon sides=5,draw,fill,inner sep=1pt,cyan] (SA2) at (axis cs: 2,80) {};
	\node [regular polygon,regular polygon sides=5,draw,fill,inner sep=1pt,cyan] (SA3) at (axis cs: 3,190) {};

	\draw [cyan] (SA1) -- (SA2) -- (SA3);
	
	\node [draw,fill,inner sep=1.25pt,blue] (DFJLS1) at (axis cs: 1,10) {};	
	\node [draw,fill,inner sep=1.25pt,blue] (DFJLS2) at (axis cs: 2,20) {};
	\node [draw,fill,inner sep=1.25pt,blue] (DFJLS3) at (axis cs: 3,40) {};
	\node [draw,fill,inner sep=1.25pt,blue] (DFJLS4) at (axis cs: 4,60) {};
	\node [draw,fill,inner sep=1.25pt,blue] (DFJLS5) at (axis cs: 5,130) {};
	\node [draw,fill,inner sep=1.25pt,blue] (DFJLS6) at (axis cs: 6,270) {};
	\node [draw,fill,inner sep=1.25pt,blue] (DFJLS7) at (axis cs: 7,490) {};
	\node [draw,fill,inner sep=1.25pt,blue] (DFJLS8) at (axis cs: 8,950) {};
	\node [draw,fill,inner sep=1.25pt,blue] (DFJLS9) at (axis cs: 9,5400) {};
	\node [draw,fill,inner sep=1.25pt,blue] (DFJLS10) at (axis cs: 10,7540) {};
	
	\draw [blue] (DFJLS1) -- (DFJLS2) -- (DFJLS3) -- (DFJLS4) -- (DFJLS5) -- (DFJLS6) -- (DFJLS7) -- (DFJLS8) -- (DFJLS9) -- (DFJLS10);
	
	\node [diamond,draw,fill,inner sep=1pt,red] (DFJ1) at (axis cs: 1,10) {};
	\node [diamond,draw,fill,inner sep=1pt,red] (DFJ2) at (axis cs: 2,70) {};
	\node [diamond,draw,fill,inner sep=1pt,red] (DFJ3) at (axis cs: 3,70) {};
	\node [diamond,draw,fill,inner sep=1pt,red] (DFJ4) at (axis cs: 4,90) {};
	\node [diamond,draw,fill,inner sep=1pt,red] (DFJ5) at (axis cs: 5,150) {};
	\node [diamond,draw,fill,inner sep=1pt,red] (DFJ6) at (axis cs: 6,180) {};
	\node [diamond,draw,fill,inner sep=1pt,red] (DFJ7) at (axis cs: 7,450) {};
	\node [diamond,draw,fill,inner sep=1pt,red] (DFJ8) at (axis cs: 8,910) {};
	\node [diamond,draw,fill,inner sep=1pt,red] (DFJ9) at (axis cs: 9,2800) {};
	\node [diamond,draw,fill,inner sep=1pt,red] (DFJ10) at (axis cs: 10,5870) {};
	
	\draw [red] (DFJ1) -- (DFJ2) -- (DFJ3) -- (DFJ4) -- (DFJ5) -- (DFJ6) -- (DFJ7) -- (DFJ8) -- (DFJ9) -- (DFJ10);

	\draw [red] (axis cs: 6.5,64) -- (axis cs: 8.5,64);
	\node [diamond,draw,fill,inner sep=1pt,red] at (axis cs: 7.5,64) {};
	\node [red] at (axis cs: 9.65,64) {\footnotesize DFJ};

	\draw [blue] (axis cs: 6.5,32) -- (axis cs: 8.5,32);
	\node [draw,fill,inner sep=1.25pt,blue] at (axis cs: 7.5,32) {};
	\node [blue] at (axis cs: 9.65,32) {\footnotesize DFJ+LS};
	
	\draw [olive] (axis cs: 6.5,16) -- (axis cs: 8.5,16);
	\node [star,star points=6,draw,fill,inner sep=1pt,olive] at (axis cs: 7.5,16) {};
	\node [olive] at (axis cs: 9.65,16) {\footnotesize MTZ};
	
	\draw [cyan] (axis cs: 6.5,8) -- (axis cs: 8.5,8);
	\node [regular polygon,regular polygon sides=5,draw,fill,inner sep=1pt,cyan] at (axis cs: 7.5,8) {};
	\node [cyan] at (axis cs: 9.65,8) {\footnotesize SA};

	\end{axis}
	\end{tikzpicture}
	\caption{Computational results for directed pyramidal tours}
	\label{image:results_directed_pyramidal}
\end{figure}
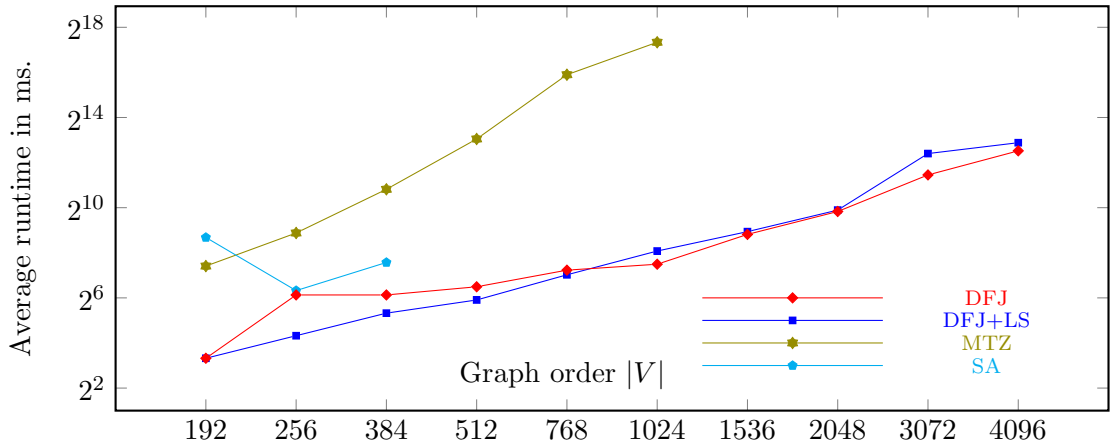

On random permutations (Fig.~\ref{image:results_directed_permutations}, Table~\ref{table:directed}, line 1) algorithms DFJ, DFJ+LS and MTZ correctly solved all $1\,000 $ instances, finding $194$ answers. The winner by margin was the Dantzig-Fulkerson-Johnson (DFJ) formulation, which found a solution in $6.15 \pm 4.3$ iterations and proved that there was no solution in $6.44 \pm 3.5$ iterations. It can be seen, that random directed multigraphs contain few subtours, which allows to quickly find a solution or prove that a solution does not exist. DFJ+LS lagged an average of $11.5 \pm 12.3$ times. On infeasible problems, local search only slows down the algorithm. But the loss on feasible problems was also significant. The third result was for the Miller-Tucker-Zemlin formulation, which lost to DFJ+LS by an average of $1.51 \pm 0.7$ times. The SA algorithm showed the last result, finding only 91 solutions out of 194.

The situation is fundamentally different on four-peak cycles (Fig.~\ref{image:results_directed_4_peaks}, Table~\ref{table:directed}, line 2). Here the winner was the DFJ+LS algorithm with the local search heuristic, which solved 897 instances out of $1\,000$. The pure DFJ lagged on average $2.4 \pm 1.7$ times and solved $798$ problems. The third result was shown by MTZ, which lagged behind DFJ by an average of $2.8 \pm 1.5$ times and solved 764 instances. The SA heuristic was able to find only 95 solutions out of 934, behind MTZ by $18.4 \pm 1.1$ times.

On pyramidal tours (Fig.~\ref{image:results_directed_pyramidal}, Table~\ref{table:directed}, line 3) DFJ and DFJ+LS solved all $1\,000$ instances out of $1\,000$, showing statistically indistinguishable results. The third place is taken by MTZ, which solved $528$ problems and lagged by $217.8 \pm 235.2$ times on average. The SA heuristic was formally faster than MTZ on feasible problems, but it was able to solve only $113$ instances out of $1\,000$.

In general, on directed graphs, based on the computational experiments, the winner was the pure Dantzig-Fulkerson-Johnson formulation, which in most cases found a solution or proved that a solution does not exist in just a few iterations.
Note that the running time and the number of iterations for feasible and infeasible problems were almost the same.
Problems arose only on four-peak cycles, on which DFJ+LS came out ahead with an additional local search heuristic. The Miller-Tucker-Zemlin model on directed graphs performed much better than on undirected ones. This was predictable since the model was originally developed for the asymmetric travelling salesman problem \cite{Miller1960}. However, it still lost significantly to the Dantzig-Fulkerson-Johnson formulation. The worst results were shown by the SA heuristic based on random matchings, which for infeasible problems could not prove that there is no solution, and for many feasible problems, it could not find a solution in the allotted time.

Some final thoughts on the computational results.
Undirected multigraphs usually contain a large number of subtours.
On the one hand, this creates significant problems for the Dantzig-Fulkerson-Johnson formulation, which tries to forbid every subtour.
On the other hand, the same fact helps heuristics that work with edge-disjoint 2-factors.
The four-peak cycles were the trickiest during testing.
First, they generate significantly fewer Hamiltonian decompositions than random permutations. Second, the small number of multiple edges does not allow us to simplify the ILP-model as much as in the case of pyramidal tours.
As for directed multigraphs, they generally do not contain many subtours, which allows the Dantzig-Fulkerson-Johnson formulation to quickly solve the problem in just a few iterations.

\section{Conclusion}

In this paper, we consider the problem of finding the second Hamiltonian decomposition of a 4-regular multigraph, motivated by verifying non-adjacency in the 1-skeleton of a travelling salesperson polytope. We adapted the classical Dantzig-Fulkerson-Johnson \cite{Dantzig1954} and Miller-Tucker-Zemlin \cite {Miller1960} formulations for the travelling salesperson problem, as well as the variable neighbourhood descent heuristic from \cite {Nikolaev2021}.

Compared to the general variable neighbourhood search algorithm from \cite {Nikolaev2021}, the new algorithm is exact and based on an integer linear programming model. It also contains an additional chain edge fixing procedure that improves performance on multigraphs with a large number of multiple edges.
Compared to the previous version of the algorithm in the proceedings of the conference ``MOTOR 2021'' \cite{Kostenko2021}, the local search has been replaced by a variable neighbourhood descent w.r.t. two neighbourhood structures.

Based on the results of computational experiments on graphs of different types, we can conclude that the integer linear programming approach works better than random matchings \cite{Nikolaev2019}. On infeasible problems, the ILP can prove that a solution does not exist. And on feasible problems, especially on undirected graphs, performance can be significantly improved by the additional heuristics.

Note that the algorithms described here can be easily adapted for the basic problem of constructing the Hamiltonian decomposition of a 4-regular multigraph and its various practical applications. It is enough to exclude from the model the constraints (\ref{DFJ_not_x})--(\ref{DFJ_not_y}), which forbid the given Hamiltonian cycles.

\section*{Acknowledgement(s)}
We acknowledge the contribution of Andrey N. Kostenko who participated in the development of the previous version of the algorithm presented at the conference ``MOTOR 2021'' \cite{Kostenko2021}.

\section*{Disclosure statement}
No potential conflict of interest was reported by the authors.

\section*{Funding}

Andrei V. Nikolaev's work was partially supported by the P.G. Demidov Yaroslavl State University Project VIP-016.

\section*{Notes on contributors}

Andrei V. Nikolaev received PhD in discrete mathematics and mathematical cybernetics from the P.G. Demidov Yaroslavl State University in 2011 and is currently an associate professor of the Department of Discrete Analysis, P.G. Demidov Yaroslavl State University. His research interest includes combinatorial optimization and polyhedral combinatorics.

Egor V. Klimov is an undergraduate student of the Faculty of Information and Computer Science, P.G. Demidov Yaroslavl State University, studying Applied Mathematics and Informatics. This work is part of his bachelor's thesis.

\section*{ORCID}

Andrei V. Nikolaev \url{https://orcid.org/0000-0003-4705-2409}\\
Egor. V. Klimov \url{https://orcid.org/0000-0002-9250-5115}

\bibliographystyle{tfs}
\bibliography{KlimovNikolaev2022}

\begin{thebibliography}{10}
\providecommand{\MR}{\relax\unskip\space MR }
\providecommand{\url}[1]{\normalfont{#1}}
\providecommand{\urlprefix}{Available at }

\bibitem{Aguilera2017}
N.E. Aguilera, R.D. Katz, and P.B. Tolomei, \emph{Vertex adjacencies in the set
  covering polyhedron}, Discrete Applied Mathematics 218 (2017), pp. 40--56.

\bibitem{Alfakih1998}
A.Y. Alfakih and K.G. Murty, \emph{Adjacency on the constrained assignment
  problem}, Discrete Applied Mathematics 87 (1998), pp. 269--274.

\bibitem{Alspach1990}
B. Alspach, J.C. Bermond, and D. Sotteau, \emph{Decomposition into Cycles I:
  Hamilton Decompositions}, in \emph{Cycles and Rays}, G. Hahn, G. Sabidussi,
  and R.E. Woodrow, eds., Springer Netherlands, Dordrecht (1990), pp. 9--18.

\bibitem{Alspach2008}
B. Alspach, \emph{The wonderful walecki construction}, Bulletin of the
  Institute of Combinatorics and its Applications 52 (2008), pp. 7--20.

\bibitem{Applegate2006}
D.L. Applegate, R.E. Bixby, V. Chvatál, and W.J. Cook, \emph{The Traveling
  Salesman Problem: A Computational Study}, Princeton University Press, 2006.

\bibitem{Arthanari2006}
T.S. Arthanari, \emph{On pedigree polytopes and hamiltonian cycles}, Discrete
  Mathematics 306 (2006), pp. 1474--1492.

\bibitem{Arthanari2013}
T.S. Arthanari, \emph{Study of the pedigree polytope and a sufficiency
  condition for nonadjacency in the tour polytope}, Discrete Optimization 10
  (2013), pp. 224--232.

\bibitem{Bae2003}
M. Bae and B. Bose, \emph{Edge disjoint hamiltonian cycles in k-ary n-cubes and
  hypercubes}, IEEE Transactions on Computers 52 (2003), pp. 1271--1284.

\bibitem{Bailey2009}
R.F. Bailey, \emph{Error-correcting codes from permutation groups}, Discrete
  Mathematics 309 (2009), pp. 4253--4265.

\bibitem{Balas1975}
E. Balas and M. Padberg, \emph{On the set-covering problem: Ii. an algorithm
  for set partitioning}, Operations Research 23 (1975), pp. 74--90.

\bibitem{Balinski1985}
M.L. Balinski, \emph{Signature methods for the assignment problem}, Operations
  Research 33 (1985), pp. 527--536.

\bibitem{Bondarenko1983}
V.A. Bondarenko, \emph{Nonpolynomial lower bounds for the complexity of the
  traveling salesman problem in a class of algorithms}, Automation and Remote
  Control 44 (1983), pp. 1137--1142.

\bibitem{Bondarenko2013}
V.A. Bondarenko and A.V. Nikolaev, \emph{Combinatorial and geometric properties
  of the max-cut and min-cut problems}, Doklady Mathematics 88 (2013), pp.
  516--517.

\bibitem{Bondarenko2018}
V.A. Bondarenko and A.V. Nikolaev, \emph{On the skeleton of the polytope of
  pyramidal tours}, Journal of Applied and Industrial Mathematics 12 (2018),
  pp. 9--18.

\bibitem{Clifton2002}
C. Clifton, M. Kantarcioglu, J. Vaidya, X. Lin, and M.Y. Zhu, \emph{Tools for
  privacy preserving distributed data mining}, SIGKDD Explor. Newsl. 4 (2002),
  pp. 28--34.

\bibitem{Cygan2015}
M. Cygan, F.V. Fomin, {\L}. Kowalik, D. Lokshtanov, D. Marx, M. Pilipczuk, M.
  Pilipczuk, and S. Saurabh, \emph{Bounded Search Trees}, in
  \emph{Parameterized Algorithms}, Springer International Publishing, Cham
  (2015), pp. 51--76.

\bibitem{Dantzig1954}
G. Dantzig, R. Fulkerson, and S. Johnson, \emph{Solution of a large-scale
  traveling-salesman problem}, Journal of the Operations Research Society of
  America 2 (1954), pp. 393--410.

\bibitem{Dantzig1963}
G.B. Dantzig, \emph{Linear Programming and Extensions}, RAND Corporation, Santa
  Monica, CA, 1963.

\bibitem{DeKort1993}
J.B. {De Kort}, \emph{A branch and bound algorithm for symmetric 2-peripatetic
  salesman problems}, European Journal of Operational Research 70 (1993), pp.
  229--243.

\bibitem{Dong2010}
R. Dong and R. Kresman, \emph{Notes on Privacy-preserving Distributed Mining
  and Hamiltonian Cycles}, in \emph{{ICSOFT} 2010 - Proceedings of the Fifth
  International Conference on Software and Data Technologies, Volume 1, Athens,
  Greece, July 22-24, 2010}, J.A.M. Cordeiro, M. Virvou, and B. Shishkov, eds.
  SciTePress, 2010, pp. 103--107.

\bibitem{Duarte2018}
A. Duarte, J. S{\'a}nchez-Oro, N. Mladenovi{\'{c}}, and R. Todosijevi{\'{c}},
  \emph{Variable Neighborhood Descent}, in \emph{Handbook of Heuristics}, R.
  Mart{\'i}, P.M. Pardalos, and M.G.C. Resende, eds., Springer International
  Publishing, Cham (2018), pp. 341--367.

\bibitem{Edmonds1965}
J. Edmonds, \emph{Paths, trees, and flowers}, Canadian Journal of Mathematics
  17 (1965), pp. 449--467.

\bibitem{Gilmore1985}
P. Gilmore, E. Lawler, and D. Shmoys, \emph{Well-solved special cases}, in
  \emph{The Traveling Salesman Problem: A Guided Tour of Combinatorial
  Optimization}, E. Lawler, J. Lenstra, A. Rinnooy~Kan, and D. Shmoys, eds.,
  Wiley (1985), pp. 87--143.

\bibitem{Glebov2017}
R. Glebov, Z. Luria, and B. Sudakov, \emph{The number of hamiltonian
  decompositions of regular graphs}, Israel Journal of Mathematics 222 (2017),
  pp. 91--108.

\bibitem{Grotchel1985}
M. Gr\"{o}tschel and M. Padberg, \emph{Polyhedral theory}, in \emph{The
  Traveling Salesman Problem: A Guided Tour of Combinatorial Optimization},
  John Wiley, Chichester,  1985, pp. 251--305.

\bibitem{Gurobi}
 {Gurobi Optimization, LLC}, \emph{Gurobi 9.1.2} (2022).
  \urlprefix\url{https://www.gurobi.com}.

\bibitem{Hansen2019}
P. Hansen, N. Mladenovi{\'{c}}, J. Brimberg, and J.A.M. P{\'e}rez,
  \emph{Variable Neighborhood Search}, in \emph{Handbook of Metaheuristics}, M.
  Gendreau and J.Y. Potvin, eds., Springer International Publishing, Cham
  (2019), pp. 57--97.

\bibitem{Hansen2017}
P. Hansen, N. Mladenovi\'c, R. Todosijevi\'c, and S. Hanafi, \emph{Variable
  neighborhood search: basics and variants}, EURO Journal on Computational
  Optimization 5 (2017), pp. 423--454.

\bibitem{Hausmann1978}
D. Hausmann and B. Korte, \emph{Colouring criteria for adjacency on
  0--1-polyhedra}, in \emph{Polyhedral Combinatorics: Dedicated to the memory
  of D.R. Fulkerson}, M.L. Balinski and A.J. Hoffman, eds., Springer Berlin
  Heidelberg, Berlin, Heidelberg (1978), pp. 106--127.

\bibitem{Hung2011}
R.W. Hung, \emph{Embedding two edge-disjoint hamiltonian cycles into locally
  twisted cubes}, Theoretical Computer Science 412 (2011), pp. 4747--4753.

\bibitem{Ikura1985}
Y. Ikura and G.L. Nemhauser, \emph{Simplex pivots on the set packing polytope},
  Mathematical Programming 33 (1985), pp. 123--138.

\bibitem{Karp1972}
R.M. Karp, \emph{Reducibility among Combinatorial Problems}, in
  \emph{Complexity of Computer Computations: Proceedings of a symposium on the
  Complexity of Computer Computations, March 20--22, 1972, New York}, R.E.
  Miller, J.W. Thatcher, and J.D. Bohlinger, eds., Springer US, Boston, MA
  (1972), pp. 85--103.

\bibitem{Kim2001}
J.H. Kim and N.C. Wormald, \emph{Random matchings which induce hamilton cycles
  and hamiltonian decompositions of random regular graphs}, Journal of
  Combinatorial Theory, Series B 81 (2001), pp. 20--44.

\bibitem{Knuth1997}
D.E. Knuth, \emph{The Art of Computer Programming, Volume 2 (3rd Ed.):
  Seminumerical Algorithms}, Addison-Wesley Longman Publishing Co., Inc., USA,
  1997.

\bibitem{Kostenko2021}
A. Kostenko and A. Nikolaev, \emph{An Iterative ILP Approach for Constructing a
  Hamiltonian Decomposition of a Regular Multigraph}, in \emph{Mathematical
  Optimization Theory and Operations Research: Recent Trends}, A. Strekalovsky,
  Y. Kochetov, T. Gruzdeva, and A. Orlov, eds., Cham. Springer International
  Publishing, 2021, pp. 216--232.

\bibitem{Kozlova2019}
A. Kozlova and A. Nikolaev, \emph{Simulated annealing approach to verify vertex
  adjacencies in the traveling salesperson polytope}, in \emph{Mathematical
  Optimization Theory and Operations Research. MOTOR 2019}, LNCS Vol. 11548.
  Springer, 2019, pp. 374--389.

\bibitem{Krarup1995}
J. Krarup, \emph{The Peripatetic Salesman and Some Related Unsolved Problems},
  in \emph{Combinatorial Programming: Methods and Applications}, B. Roy, ed.,
  Vol.~19, Dordrecht. Springer Netherlands, 1995, pp. 173--178.

\bibitem{Michini2014}
C. Michini and A. Sassano, \emph{The hirsch conjecture for the fractional
  stable set polytope}, Mathematical Programming 147 (2014), pp. 309--330.

\bibitem{Miller1960}
C.E. Miller, A.W. Tucker, and R.A. Zemlin, \emph{Integer programming
  formulation of traveling salesman problems}, J. ACM 7 (1960), pp. 326--329.

\bibitem{Mladenovic1997}
N. Mladenovi\'c and P. Hansen, \emph{Variable neighborhood search}, Computers
  \& Operations Research 24 (1997), pp. 1097--1100.

\bibitem{Nikolaev2019}
A. Nikolaev, \emph{On vertex adjacencies in the polytope of pyramidal tours
  with step-backs}, in \emph{Mathematical Optimization Theory and Operations
  Research. MOTOR 2019}, LNCS Vol. 11548. Springer, 2019, pp. 247--263.

\bibitem{Nikolaev2021}
A. Nikolaev and A. Kozlova, \emph{Hamiltonian decomposition and verifying
  vertex adjacency in 1-skeleton of the traveling salesperson polytope by
  variable neighborhood search}, Journal of Combinatorial Optimization 42
  (2021), pp. 212--230.

\bibitem{Orman2007}
A. Orman and H. Williams, \emph{A Survey of Different Integer Programming
  Formulations of the Travelling Salesman Problem}, in \emph{Optimisation,
  Econometric and Financial Analysis}, E.J. Kontoghiorghes and C. Gatu, eds.,
  Berlin, Heidelberg. Springer Berlin Heidelberg, 2007, pp. 91--104.

\bibitem{Papadimitriou1978}
C.H. Papadimitriou, \emph{The adjacency relation on the traveling salesman
  polytope is np-complete}, Mathematical Programming 14 (1978), pp. 312--324.

\bibitem{Peroche1984}
B. P\'{e}roche, \emph{Np-completeness of some problems of partitioning and
  covering in graphs}, Discrete Applied Mathematics 8 (1984), pp. 195--208.

\bibitem{Rao1976}
M.R. Rao, \emph{Adjacency of the traveling salesman tours and \$0 - 1\$
  vertices}, SIAM Journal on Applied Mathematics 30 (1976), pp. 191--198.

\bibitem{Rowley1991}
R. Rowley and B. Bose, \emph{Edge-disjoint hamiltonian cycles in de bruijn
  networks} (1991).

\bibitem{Tutte1954}
W.T. Tutte, \emph{A short proof of the factor theorem for finite graphs},
  Canadian Journal of Mathematics 6 (1954), pp. 347--352.

\end{thebibliography}

\end{document}